\documentclass[11pt,reqno]{amsart}
\usepackage[utf8]{inputenc}
\usepackage[T1]{fontenc}
\usepackage[margin=1in]{geometry}
\usepackage{amsmath,amssymb,amsthm,mathtools}
\usepackage{graphicx,xcolor}
\usepackage{circuitikz}
\usepackage{tikz-3dplot}
\usepackage{subcaption}
\usepackage{pgfplots}
\usepackage{float}
\usepackage[colorlinks=true,linkcolor=black,citecolor=black,urlcolor=blue]{hyperref}
\usetikzlibrary{calc,positioning,lindenmayersystems}
\pgfplotsset{compat=1.18}
\numberwithin{equation}{section}
\newtheorem{theorem}{Theorem}[section]
\newtheorem{proposition}[theorem]{Proposition}
\newtheorem{lemma}[theorem]{Lemma}
\newtheorem{corollary}[theorem]{Corollary}
\theoremstyle{definition}
\newtheorem{definition}[theorem]{Definition}

\theoremstyle{remark}

\allowdisplaybreaks
\title{The Fourier Ratio on Resistance Spaces and Applications to Electrical Networks}
\author{Christopher Housholder}
\address{Department of Mathematics, Missouri State University, Springfield, MO 65897}
\email{christopherlancehousholder@gmail.com}
\author{Alex Iosevich}
\address{Department of Mathematics, University of Rochester, Rochester, NY 14627}
\email{alex.iosevich@rochester.edu}
\thanks{The second author was supported in part by NSF grant DMS-2506858 and a Simons Fellowship in Mathematics in 2026.}
\author{Steven J. Miller}
\address{Department of Mathematics, Williams College, Williamstown, MA 01267}
\email{sjm1@williams.edu}
\author{Eyvindur A. Palsson}
\address{Department of Mathematics, Virginia Tech, Blacksburg, VA 24061}
\email{palsson@vt.edu}
\date{October 19, 2026}
\begin{document}
\begin{abstract}
We study the Fourier Ratio for spectral expansions of nonnegative self-adjoint operators. Weighted Cauchy--Schwarz estimates turn eigenvalue counting bounds into Fourier Ratio bounds for spectrally regular functions. On a compact measured resistance space, Kasue's identity expresses the reciprocal eigenvalue sum as one half of the average resistance distance. We use this identity to obtain a sharp coefficient estimate for every finite-energy function, with no spectral cutoff. Bounded orthonormal system recovery theorems then yield reconstruction estimates whose sample counts depend on both the Fourier Ratio and the pointwise size of the chosen eigenfunctions. We distinguish these estimates from the weaker diameter-only bounds and give an interval example with a sample count logarithmic in the dimension of the spectral subspace for fixed accuracy and energy ratio. Resistance also controls errors caused by displaced sensors. We apply the estimates to standard self-similar fractals and finite electrical networks, where the geometric constant becomes the normalized square root of the Kirchhoff index. A separate observation describes Fourier Ratios of coefficient families with a limiting distribution.
\end{abstract}
\maketitle
\tableofcontents

\section{Introduction}
The Fourier Ratio measures the concentration of a coefficient vector. For a nonzero function with Fourier coefficients $\widehat{f}$, it is defined by
\begin{equation}
\operatorname{FR}(f)=\frac{\left\|\widehat{f}\right\|_1}{\left\|\widehat{f}\right\|_2}.
\end{equation}
If exactly $s$ coefficients have the same nonzero magnitude, then the ratio is $\sqrt{s}$. For a general vector, it measures spectral concentration without requiring exact sparsity. Recent work has used this quantity in signal recovery, uncertainty principles, restriction estimates, sampling, and learning theory \cite{AldalehEtAl,BursteinIosevichNathan,IosevichLiPalssonYavicoli,IosevichPalssonYavicoli,IosevichIosevichPalssonYavicoli}. These connections belong to the broader theory of uncertainty and recovery from incomplete measurements \cite{DonohoStark,CandesRombergTaoRobust,CandesRombergTaoStable,FoucartRauhut}.

Our starting question is how geometry and regularity can control the Fourier Ratio. Suppose that $A$ is a nonnegative self-adjoint operator with discrete spectrum and that $f=\sum_j\widehat{f}(j)e_j$ in a fixed orthonormal eigenbasis. After separating the constant eigenfunction, weighted Cauchy--Schwarz gives
\begin{equation}\label{eq:introspectral}
\sum_{j\geq1}\left|\widehat{f}(j)\right|
\leq
\left(\sum_{j\geq1}\lambda_j^{-s}\right)^{1/2}
\left\|A^{s/2}f\right\|_2,
\end{equation}
whenever the right side is finite. Thus regularity controls the coefficients once the reciprocal spectral sum is understood. Weyl-type counting estimates supply one way to estimate this sum \cite{Weyl1911,HormanderSpectralFunction,SafarovVassiliev}. The three possibilities are familiar. The sum is bounded, logarithmic, or grows like a power, according to the relation between the regularity order and the counting exponent.

On compact measured resistance spaces, the spectral sum has a direct geometric interpretation. A resistance form supplies an energy $\mathcal E$, a domain $\mathcal F$, and a resistance metric $R$ \cite{KigamiResistance,KigamiAnalysis}. With a full-support probability measure $\mu$, the associated operator has constants as its kernel. Kasue proved the identities
\begin{equation}\label{eq:introtrace}
R(x,y)=\sum_{j\geq1}\frac{\left|e_j(x)-e_j(y)\right|^2}{\lambda_j},
\qquad
\sum_{j\geq1}\frac{1}{\lambda_j}
=\frac12\int_X\int_X R(x,y)\,d\mu(x)\,d\mu(y).
\end{equation}
See \cite[the equations (3.3) and (3.5)]{Kasue}. We include a proof in our normalization because the factor $1/2$ enters the Fourier Ratio and network estimates.

Write
\begin{equation}\label{eq:introCR}
C_R=\left(\frac12\int_X\int_XR(x,y)\,d\mu(x)\,d\mu(y)\right)^{1/2}.
\end{equation}
Our principal use of the identity is the sharp estimate
\begin{equation}\label{eq:introsharp}
\sum_{j\geq1}\left|\widehat{f}(j)\right|
\leq C_R\mathcal E(f,f)^{1/2},\qquad f\in\mathcal F.
\end{equation}
The constant is optimal for every fixed orthonormal eigenbasis. In particular, every finite-energy function has an absolutely summable spectral coefficient sequence, and
\begin{equation}
\operatorname{FR}(f)
\leq
\frac{\left|\int_Xf\,d\mu\right|}{\left\|f\right\|_2}
+C_R\frac{\mathcal E(f,f)^{1/2}}{\left\|f\right\|_2}.
\end{equation}
There is no cutoff in this statement. The dependence on the function is through its mean and its energy relative to its $L^2$ norm.

A Fourier Ratio bound is one ingredient in a point-sampling theorem. The other is the size of the chosen eigenfunctions. For the spectral space
\begin{equation}
V_L=\operatorname{span}\{e_j:\lambda_j\leq L\},
\end{equation}
let $D_L=\dim V_L$ and
\begin{equation}
K_L=\max_{\lambda_j\leq L}\left\|e_j\right\|_\infty.
\end{equation}
Standard recovery estimates for bounded orthonormal systems \cite{BrugiapagliaDirksenJungRauhut,FoucartRauhut,RauhutBOS,RauhutWard} combine these two parameters. A function with Fourier Ratio at most $r$ can be approximated to relative accuracy $\varepsilon$ from a number of independent point samples on the scale
\begin{equation}\label{eq:introsamples}
K_L^2\frac{r^2}{\varepsilon^2}
\log^2\left(2+K_L^2\frac{r^2}{\varepsilon^2}\right)\log(2D_L).
\end{equation}
The reconstruction uses $\ell^1$ minimization and is stable under measurement noise.

The distinction between the two parameters is essential. Resistance geometry gives $K_L^2\leq\max\{1,D_RL\}$, where $D_R$ is the resistance diameter. However, the same trace identity implies $D_L-1\leq D_RL/2$. The resulting diameter-only sample count therefore does not establish a saving over the number of unknown coefficients. A sharper estimate for the actual $K_L$ is needed for that purpose. On the Neumann interval, for example, $K_L^2=2$. For fixed accuracy and a bounded energy-to-norm ratio, the sample count in \eqref{eq:introsamples} is then logarithmic in $D_L$.

Fractals provide examples in which both the resistance structure and eigenvalue growth are well understood \cite{FukushimaShima,KigamiLapidus,KumagaiSpectralDimension,StrichartzSpectral,StrichartzAnalysis,Riane}. The Sierpinski gasket has counting exponent $\log 3/\log 5$, the Sierpinski tetrahedron has exponent $\log 4/\log 6$, and the Minkowski curve has exponent $1/2$. Their Fourier Ratio estimates follow directly from \eqref{eq:introsharp}. Their reconstruction estimates retain the dependence on eigenfunction size. We do not infer compressed sampling from spectral dimension alone.

Finite electrical networks give another interpretation. The resistance metric becomes effective electrical resistance, and the energy is the Joule power dissipated by a voltage potential \cite{DoyleSnell,KleinRandic,LyonsPeres}. With the uniform probability measure on $n$ vertices, the resistance-space operator is $n$ times the physical graph Laplacian. In this normalization, $C_R=\sqrt{Kf(G)}/n$, where $Kf(G)$ is the Kirchhoff index. The current-voltage equation adds an inverse-Laplacian relation, from which we obtain spectral approximation estimates for voltages that are not exactly bandlimited. The resistance inequality also bounds the change in a voltage reading when a sensor is moved.

We first establish the resistance-space setting and the spectral coefficient estimates. We then treat sampling and sensor stability, followed by fractal examples and electrical networks. The appendix records the elementary distributional observation from the original finite-family formulation.

\section{Resistance-space tools}
A resistance form is an energy for which differences of point evaluations are continuous in the energy norm. This is the feature needed below. We recall the definition and the compactness argument so that the operator assumptions do not remain implicit.

\begin{definition}\label{def:resform}
Let $X$ be a set. A pair $(\mathcal E,\mathcal F)$ is a resistance form on $X$ if the following properties hold.
\begin{enumerate}
\item $\mathcal F$ is a linear space of real-valued functions containing the constants. The form $\mathcal E$ is nonnegative, symmetric, and bilinear, and $\mathcal E(u,u)=0$ if and only if $u$ is constant.
\item After functions differing by a constant are identified, $\mathcal F$ is a Hilbert space with inner product $\mathcal E$.
\item For every finite set $V\subset X$ and every $v:V\rightarrow\mathbb{R}$, there is $u\in\mathcal F$ with $u|_V=v$.
\item For every $x,y\in X$,
\begin{equation}
R(x,y)=\sup_{\substack{u\in\mathcal F\\\mathcal E(u,u)>0}}
\frac{\left|u(x)-u(y)\right|^2}{\mathcal E(u,u)}<\infty.
\end{equation}
\item If $u\in\mathcal F$ and $\overline u=\max\{0,\min\{1,u\}\}$, then $\overline u\in\mathcal F$ and $\mathcal E(\overline u,\overline u)\leq\mathcal E(u,u)$.
\end{enumerate}
The function $R$ is a metric, called the resistance metric \cite{KigamiResistance,KigamiAnalysis}.
\end{definition}

The definition immediately implies the resistance inequality
\begin{equation}\label{eq:resistanceinequality}
\left|u(x)-u(y)\right|^2\leq R(x,y)\mathcal E(u,u).
\end{equation}
It remains valid after complexification, with the usual sesquilinear extension of $\mathcal E$. We use the same notation for this extension.

\begin{definition}\label{def:measuredresistance}
A compact measured resistance space is a tuple $(X,R,\mu,\mathcal E,\mathcal F)$ in which $(\mathcal E,\mathcal F)$ is a resistance form, $(X,R)$ is compact, and $\mu$ is a Borel probability measure with full support. Its resistance diameter is
\begin{equation}
D_R=\operatorname{diam}_R(X).
\end{equation}
\end{definition}

We assume that $X$ has at least two points whenever a first positive eigenvalue is used. For a one-point space, the energy vanishes and the Fourier Ratio of every nonzero function is $1$. Finite resistance spaces are otherwise included throughout.

\begin{lemma}\label{lem:resistanceuniform}
Let $u\in\mathcal F$ and $u_X=\int_Xu\,d\mu$. Then
\begin{equation}\label{eq:uniformmean}
\left\|u-u_X\right\|_\infty\leq D_R^{1/2}\mathcal E(u,u)^{1/2},
\end{equation}
and
\begin{equation}
\left\|u\right\|_\infty\leq\left\|u\right\|_2+D_R^{1/2}\mathcal E(u,u)^{1/2}.
\end{equation}
\end{lemma}
\begin{proof}
Since $\mu(X)=1$,
\begin{equation}
u(x)-u_X=\int_X\left(u(x)-u(y)\right)\,d\mu(y).
\end{equation}
The resistance inequality bounds the absolute value by $D_R^{1/2}\mathcal E(u,u)^{1/2}$. Taking the supremum proves \eqref{eq:uniformmean}. The second estimate follows from $\left|u_X\right|\leq\left\|u\right\|_2$.
\end{proof}

\begin{lemma}\label{lem:closedform}
The space $\mathcal F$ is dense in $L^2(X,\mu)$ and complete in the form norm
\begin{equation}
\left\|u\right\|_{\mathcal F}^2=\mathcal E(u,u)+\left\|u\right\|_2^2.
\end{equation}
\end{lemma}
\begin{proof}
Every finite-energy function is continuous by \eqref{eq:resistanceinequality}, and is therefore bounded. The Markov property implies that $u^+\in\mathcal F$ for bounded real $u\in\mathcal F$, by rescaling the truncation in Definition~\ref{def:resform}. Thus $\mathcal F$ is a vector lattice containing the constants. It separates points by the finite interpolation axiom. The lattice form of the Stone--Weierstrass theorem shows that it is uniformly dense in $C(X)$, and hence dense in $L^2(X,\mu)$.

Let $\mathcal F_0$ be the mean-zero subspace of $\mathcal F$. Choosing the mean-zero representative identifies $\mathcal F_0$ with the Hilbert space of functions modulo constants, with inner product $\mathcal E$. By \eqref{eq:uniformmean}, convergence in this energy norm implies uniform convergence of the mean-zero representatives. The form norm is equivalent on $\mathcal F_0$ to the energy norm. Adding the constant component proves completeness on $\mathcal F$.
\end{proof}

The density and form-norm completeness statements above also hold after complexification. Indeed, $\mathcal F_{\mathbb C}=\mathcal F+i\mathcal F$ is dense in $L^2(X,\mu;\mathbb C)$, and the sesquilinear extension of $\mathcal E$ is complete in the corresponding form norm, by applying the real-valued statements separately to real and imaginary parts. In the spectral results below, we use the same notation $\mathcal F$ for this complexified form domain.

\begin{theorem}\label{thm:compactembedding}
Every bounded sequence in $\mathcal F$, equipped with the form norm, has a uniformly convergent subsequence. In particular, the embedding $\mathcal F\hookrightarrow L^2(X,\mu)$ is compact.
\end{theorem}
\begin{proof}
Lemma~\ref{lem:resistanceuniform} gives uniform boundedness. If the energies are bounded by $M$, then
\begin{equation}
\left|u(x)-u(y)\right|\leq M^{1/2}R(x,y)^{1/2},
\end{equation}
so the sequence is equicontinuous in the resistance metric. Arzel\`a--Ascoli gives a uniformly convergent subsequence. Uniform convergence also implies $L^2$ convergence because $\mu(X)=1$.
\end{proof}

\section{Spectral and resistance spaces}
\subsection{Spectral Fourier Ratio bounds}
We first keep the spectral weights arbitrary. This separates the coefficient estimate from the particular regularity norm used to control it.

\begin{theorem}\label{thm:spectral}
Let $(X,\mu)$ be a probability space and let $A$ be a densely defined nonnegative self-adjoint operator on $L^2(X,\mu)$ with compact resolvent and $\ker A=\operatorname{span}\{1\}$. Fix an orthonormal eigenbasis $\{e_j\}_{j\in J}$ with $e_0=1$ and
\begin{equation}
Ae_j=\lambda_je_j,\qquad 0=\lambda_0<\lambda_1\leq\lambda_2\leq\cdots.
\end{equation}
The index set $J$ is finite if $L^2(X,\mu)$ is finite-dimensional and is $\{0,1,2,\ldots\}$ otherwise. In the latter case, $\lambda_j\rightarrow\infty$. All spectral sums are understood to run over the chosen index set $J$.
For $L\geq\lambda_1$, set
\begin{equation}
V_L=\operatorname{span}\{e_j:\lambda_j\leq L\},\qquad
\widehat{f}(j)=\langle f,e_j\rangle,
\end{equation}
and for nonzero $f\in V_L$ define
\begin{equation}
\operatorname{FR}_L(f)=\frac{\sum_{\lambda_j\leq L}\left|\widehat{f}(j)\right|}{\left\|f\right\|_2}.
\end{equation}
For every choice of positive weights $w_j$ on the positive eigenvalues below $L$,
\begin{align}
\operatorname{FR}_L(f)
\leq{}&\frac{\left|\int_Xf\,d\mu\right|}{\left\|f\right\|_2}\notag\\
&+\left(\sum_{0<\lambda_j\leq L}w_j^{-2}\right)^{1/2}
\frac{\left(\sum_{0<\lambda_j\leq L}w_j^2\left|\widehat{f}(j)\right|^2\right)^{1/2}}{\left\|f\right\|_2}.
\end{align}
In particular, for every $s>0$,
\begin{equation}\label{eq:spectralFR}
\operatorname{FR}_L(f)\leq
\frac{\left|\int_Xf\,d\mu\right|}{\left\|f\right\|_2}
+\left(\sum_{0<\lambda_j\leq L}\lambda_j^{-s}\right)^{1/2}
\frac{\left\|A^{s/2}f\right\|_2}{\left\|f\right\|_2}.
\end{equation}
\end{theorem}
\begin{proof}
The constant coefficient is $\widehat{f}(0)=\int_Xf\,d\mu$. For the remaining coefficients, Cauchy--Schwarz gives
\begin{equation}
\sum_{0<\lambda_j\leq L}\left|\widehat{f}(j)\right|
\leq\left(\sum_{0<\lambda_j\leq L}w_j^{-2}\right)^{1/2}
\left(\sum_{0<\lambda_j\leq L}w_j^2\left|\widehat{f}(j)\right|^2\right)^{1/2}.
\end{equation}
Adding the constant coefficient and dividing by $\left\|f\right\|_2$ proves the first estimate. Taking $w_j=\lambda_j^{s/2}$ identifies the last sum with $\left\|A^{s/2}f\right\|_2^2$.
\end{proof}

The estimate is relative to a fixed eigenbasis. Rotations within an eigenspace of multiplicity greater than one can change the Fourier Ratio, although they do not change the energy or the reciprocal spectral sum. Every estimate below holds for each such choice of basis. The pointwise eigenfunction bound in the sampling theorem must refer to the same choice.

\begin{theorem}\label{thm:counting}
Under the hypotheses of Theorem~\ref{thm:spectral}, let
\begin{equation}
N_A(\Lambda)=\#\{j\in J:\lambda_j\leq\Lambda\}.
\end{equation}
Suppose that $N_A(\Lambda)\leq C_N\Lambda^\alpha$ for every $\Lambda\geq\lambda_1$, where $\alpha>0$. Then for $s>0$ and $L\geq\lambda_1$,
\begin{equation}\label{eq:reciprocalcounting}
\sum_{0<\lambda_j\leq L}\lambda_j^{-s}\leq C
\begin{cases}
1,&s>\alpha,\\
1+\log(1+L/\lambda_1),&s=\alpha,\\
L^{\alpha-s},&0<s<\alpha.
\end{cases}
\end{equation}
Here $C$ depends only on $\alpha,s,C_N$, and $\lambda_1$. Consequently,
\begin{align}
\operatorname{FR}_L(f)\leq{}&\frac{\left|\int_Xf\,d\mu\right|}{\left\|f\right\|_2}\notag\\
&+C\frac{\left\|A^{s/2}f\right\|_2}{\left\|f\right\|_2}
\begin{cases}
1,&s>\alpha,\\
\left(1+\log(1+L/\lambda_1)\right)^{1/2},&s=\alpha,\\
L^{(\alpha-s)/2},&0<s<\alpha.
\end{cases}
\end{align}
\end{theorem}
\begin{proof}
Split the positive eigenvalues into intervals $[2^n\lambda_1,2^{n+1}\lambda_1)$. The contribution from the $n$th interval is at most
\begin{equation}
(2^n\lambda_1)^{-s}N_A(2^{n+1}\lambda_1)
\leq C_N2^\alpha\lambda_1^{\alpha-s}2^{n(\alpha-s)}.
\end{equation}
If $s>\alpha$, these contributions form a decreasing geometric series. If $s=\alpha$, each interval contributes at most a constant and there are $O(1+\log(1+L/\lambda_1))$ intervals. If $s<\alpha$, the last interval controls the geometric sum. This proves \eqref{eq:reciprocalcounting}. The Fourier Ratio estimate follows from \eqref{eq:spectralFR}.
\end{proof}

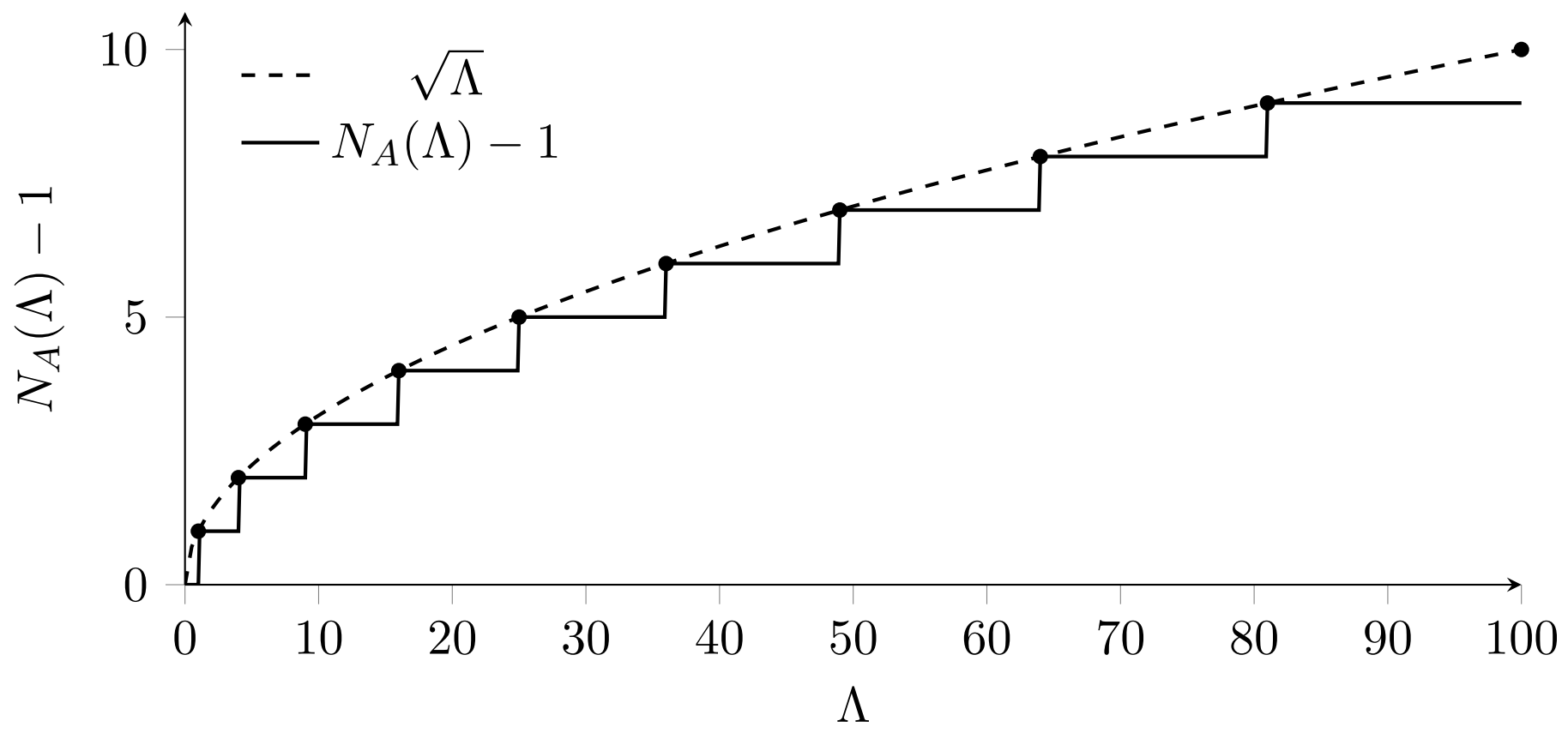
\begin{figure}[H]
\centering

\begin{tikzpicture}
\begin{axis}[
    width=0.72\textwidth,
    height=6cm,
    xlabel={$\Lambda$},
    ylabel={$N_A(\Lambda)-1$},
    xmin=0,
    xmax=100,
    ymin=0,
    ymax=10.7,
    axis lines=left,
    tick align=outside,
    legend pos=north west,
    legend style={draw=none},
    samples=200
]

% Weyl asymptotic
\addplot[
    thick,
    dashed,
    domain=0:100
]
{sqrt(x)};

\addlegendentry{$\sqrt{\Lambda}$}

% Exact counting function N_A(Lambda) = floor(sqrt(Lambda))
\addplot[
    thick,
    domain=0:100,
    samples=1001
]
{floor(sqrt(x))};

\addlegendentry{$N_A(\Lambda)-1$}

% Eigenvalue locations
\addplot[
    only marks,
    mark=*,
    mark size=1.5pt
]
coordinates {
    (1,1)
    (4,2)
    (9,3)
    (16,4)
    (25,5)
    (36,6)
    (49,7)
    (64,8)
    (81,9)
    (100,10)
};

\end{axis}
\end{tikzpicture}

\caption{The positive eigenvalues of the Neumann Laplacian on $[0,\pi]$ are $\lambda_j=j^2$. The positive counting function $N_A(\Lambda)-1=\lfloor\sqrt{\Lambda}\rfloor$ is a staircase with Weyl growth $\sqrt{\Lambda}$. The subtraction removes the constant eigenfunction.}
\label{fig:weyl-counting}
\end{figure}

\subsection{Resistance Laplacians and energy}
The operator associated with the full resistance form is the Neumann operator when the space has a distinguished boundary. No boundary values are imposed here. A Dirichlet realization, whose kernel does not contain the constants, requires a different formulation.

\begin{theorem}\label{thm:reslap}
Let $(X,R,\mu,\mathcal E,\mathcal F)$ be a compact measured resistance space. The nonnegative self-adjoint operator $A=-\Delta_\mu$ associated with $(\mathcal E,\mathcal F)$ has compact resolvent. Its kernel is the constants, and it admits the eigenbasis described in Theorem~\ref{thm:spectral}, with a finite index set when $X$ is finite. Every eigenfunction belongs to $\mathcal F\subset C(X)$, and
\begin{equation}\label{eq:formidentity}
\mathcal E(f,f)=\left\|A^{1/2}f\right\|_2^2
=\sum_{j\in J}\lambda_j\left|\widehat{f}(j)\right|^2,
\qquad f\in\mathcal F.
\end{equation}
\end{theorem}
\begin{proof}
Lemma~\ref{lem:closedform} permits the representation theorem for densely defined closed nonnegative forms \cite{Fukushima}. It supplies a self-adjoint operator satisfying
\begin{equation}
\mathcal E(u,v)=\langle Au,v\rangle,
\qquad u\in\operatorname{Dom}(A),\quad v\in\mathcal F.
\end{equation}
If $u=(A+I)^{-1}g$, then
\begin{equation}
\left\|u\right\|_{\mathcal F}^2=\langle g,u\rangle
\leq\left\|g\right\|_2\left\|u\right\|_{\mathcal F}.
\end{equation}
Thus the resolvent maps bounded subsets of $L^2$ to bounded subsets of $\mathcal F$. Theorem~\ref{thm:compactembedding} makes it compact as an operator on $L^2$. The spectral theorem now supplies the stated eigenbasis. The resistance form axioms identify the zero eigenspace with the constants. Finally, the representation theorem identifies $\mathcal F$ with $\operatorname{Dom}(A^{1/2})$ and proves \eqref{eq:formidentity}.
\end{proof}

\begin{theorem}\label{thm:energy}
Every nonzero $f\in V_L$ on a compact measured resistance space satisfies
\begin{equation}
\operatorname{FR}_L(f)\leq
\frac{\left|\int_Xf\,d\mu\right|}{\left\|f\right\|_2}
+\left(\sum_{0<\lambda_j\leq L}\frac{1}{\lambda_j}\right)^{1/2}
\frac{\mathcal E(f,f)^{1/2}}{\left\|f\right\|_2}.
\end{equation}
The first term is at most $1$.
\end{theorem}
\begin{proof}
Take $s=1$ in Theorem~\ref{thm:spectral} and use \eqref{eq:formidentity}. Cauchy--Schwarz and $\mu(X)=1$ bound the mean term by $1$.
\end{proof}

\subsection{The reciprocal spectrum and resistance}
The preceding theorem leaves a reciprocal spectral sum. Kasue's identity expresses the full sum in terms of resistance, without an eigenvalue counting assumption.

\begin{theorem}[Kasue]\label{thm:resistancetrace}
On a compact measured resistance space,
\begin{equation}\label{eq:resistanceseries}
R(x,y)=\sum_{\substack{j\in J\\j\geq1}}
\frac{\left|e_j(x)-e_j(y)\right|^2}{\lambda_j},
\qquad x,y\in X,
\end{equation}
and
\begin{equation}\label{eq:resistancetrace}
\sum_{\substack{j\in J\\j\geq1}}\frac{1}{\lambda_j}
=\frac12\int_X\int_XR(x,y)\,d\mu(x)\,d\mu(y)
\leq\frac{D_R}{2}.
\end{equation}
These are \cite[the equations (3.3) and (3.5)]{Kasue} in the probability-measure normalization.
\end{theorem}
\begin{proof}
Work in the mean-zero energy space $\mathcal F_0$. The functions
\begin{equation}
\phi_j=\frac{e_j}{\sqrt{\lambda_j}},\qquad j\geq1,
\end{equation}
form an orthonormal basis with respect to $\mathcal E$, by \eqref{eq:formidentity}. The functional $u\mapsto u(x)-u(y)$ has squared norm $R(x,y)$ in this Hilbert space. Parseval for its Riesz representative proves \eqref{eq:resistanceseries}.

All summands are nonnegative, so Tonelli's theorem permits integration over $X\times X$. Since $e_j$ has mean zero and $L^2$ norm $1$ for $j\geq1$,
\begin{equation}
\int_X\int_X\left|e_j(x)-e_j(y)\right|^2\,d\mu(x)\,d\mu(y)=2.
\end{equation}
This proves the equality in \eqref{eq:resistancetrace}. The diameter bound follows from $R(x,y)\leq D_R$.
\end{proof}

\begin{corollary}\label{cor:uniformresistance}
Let $C_R$ be defined by \eqref{eq:introCR}. Every $f\in\mathcal F$ satisfies the sharp coefficient estimate
\begin{equation}\label{eq:sharpcoefficient}
\sum_{\substack{j\in J\\j\geq1}}\left|\widehat{f}(j)\right|
\leq C_R\mathcal E(f,f)^{1/2}.
\end{equation}
Consequently, for every nonzero $f\in\mathcal F$, the full Fourier Ratio is finite and
\begin{align}\label{eq:fullresistanceFR}
\operatorname{FR}(f)
&\leq\frac{\left|\int_Xf\,d\mu\right|}{\left\|f\right\|_2}
+C_R\frac{\mathcal E(f,f)^{1/2}}{\left\|f\right\|_2}\\
&\leq1+\left(\frac{D_R}{2}\right)^{1/2}
\frac{\mathcal E(f,f)^{1/2}}{\left\|f\right\|_2}.\notag
\end{align}
The constant $C_R$ in \eqref{eq:sharpcoefficient} is optimal. For $f\in V_L$, the full ratio agrees with $\operatorname{FR}_L(f)$.
\end{corollary}
\begin{proof}
Weighted Cauchy--Schwarz on finite partial sums gives
\begin{equation}
\sum_{j\geq1}\left|\widehat{f}(j)\right|
\leq\left(\sum_{j\geq1}\lambda_j^{-1}\right)^{1/2}
\left(\sum_{j\geq1}\lambda_j\left|\widehat{f}(j)\right|^2\right)^{1/2}.
\end{equation}
Theorem~\ref{thm:resistancetrace} and \eqref{eq:formidentity} prove \eqref{eq:sharpcoefficient}. Adding the constant coefficient proves \eqref{eq:fullresistanceFR}.

For sharpness, choose the mean-zero function with coefficients $\widehat{f}(j)=1/\lambda_j$ for $j\geq1$. It belongs to $\mathcal F$ because its energy is $\sum_{j\geq1}1/\lambda_j<\infty$, and it belongs to $L^2$ because $\lambda_j\geq\lambda_1>0$. Equality holds in \eqref{eq:sharpcoefficient}.
\end{proof}

\begin{proposition}\label{prop:uniformspectralconvergence}
Every finite-energy spectral expansion converges absolutely and uniformly on $X$.
\end{proposition}
\begin{proof}
The evaluation functional $u\mapsto u(x)$ on $\mathcal F_0$ has squared norm at most $D_R$, by Lemma~\ref{lem:resistanceuniform}. Parseval in the energy basis therefore gives
\begin{equation}
\sum_{j\geq1}\frac{\left|e_j(x)\right|^2}{\lambda_j}\leq D_R.
\end{equation}
For every $M$,
\begin{equation}
\sum_{j>M}\left|\widehat{f}(j)e_j(x)\right|
\leq D_R^{1/2}
\left(\sum_{j>M}\lambda_j\left|\widehat{f}(j)\right|^2\right)^{1/2}.
\end{equation}
The right side tends to zero independently of $x$. The uniform limit agrees with $f$ in $L^2$, and hence everywhere by continuity and full support of $\mu$.
\end{proof}

\begin{corollary}\label{cor:sublinearcounting}
The counting function of a compact measured resistance space satisfies
\begin{equation}
N_A(\Lambda)=o(\Lambda),\qquad\Lambda\rightarrow\infty.
\end{equation}
\end{corollary}
\begin{proof}
The assertion is immediate for a finite space. Otherwise, given $\varepsilon>0$, choose $M$ so that $\sum_{j>M}1/\lambda_j<\varepsilon$. Then
\begin{equation}
\frac{N_A(\Lambda)-1-M}{\Lambda}
\leq\sum_{\substack{j>M\\\lambda_j\leq\Lambda}}\frac1{\lambda_j}
<\varepsilon
\end{equation}
whenever $\Lambda>\lambda_M$. Letting $\Lambda$ tend to infinity and then $\varepsilon$ tend to zero proves the claim.
\end{proof}

The absence of a cutoff in \eqref{eq:fullresistanceFR} concerns the constant multiplying the energy ratio. It is not a uniform bound over all unit vectors in $V_L$, whose energies can grow with $L$.

% ========================= Sampling and Reconstruction =========================
\section{Sampling and Reconstruction}
The spectral estimates above control the Fourier Ratio through regularity and resistance geometry. To use them for reconstruction, we also need to control the size of the eigenfunctions at the sampling locations. These are separate requirements. A small Fourier Ratio controls the concentration of the coefficient vector, while a bound on the eigenfunctions controls the interaction between point evaluation and that vector.

We begin with the standard recovery theorem for a bounded orthonormal system. We state the result with measurement noise so that the same theorem can be used for perturbed sensors and for voltages that are only approximately confined to a spectral subspace. The sampling estimate follows from \cite[Theorems~1.1 and~2.3]{BrugiapagliaDirksenJungRauhut} and stable basis pursuit \cite{FoucartRauhut}. The Fourier Ratio enters through the usual best-term approximation estimate, as in \cite{BursteinIosevichNathan}.

% ========================= Spectral Recovery =========================
\subsection{Spectral Recovery}
\begin{theorem}\label{thm:spectralrecovery}
Assume the hypotheses of Theorem~\ref{thm:spectral}. For a spectral cutoff $L$, define
\begin{equation}
D_L \ = \ \#\{j:\lambda_j\leq L\}, \qquad K_L \ = \ \max_{\lambda_j\leq L}\left\|e_j\right\|_{L^\infty(X)},
\end{equation}
and assume $K_L<\infty$. Let $0<\varepsilon<1/2$, let $f\in V_L$ be nonzero, and suppose that
\begin{equation}
\operatorname{FR}_L(f) \ \leq \ r.
\end{equation}
Choose $x_1,\ldots,x_q$ independently according to $\mu$. Suppose that the observations are
\begin{equation}
g_i \ = \ f(x_i)+\xi_i, \qquad 1\leq i\leq q,
\end{equation}
and that a prescribed radius $\eta\geq0$ satisfies
\begin{equation}
\left(\frac1q\sum_{i=1}^q\left|\xi_i\right|^2\right)^{1/2} \ \leq \ \eta.
\end{equation}
Let $f^*\in V_L$ minimize
\begin{equation}
\sum_{\lambda_j\leq L}\left|\widehat{h}(j)\right|
\end{equation}
over $h\in V_L$, subject to
\begin{equation}
\left(\frac1q\sum_{i=1}^q\left|h(x_i)-g_i\right|^2\right)^{1/2} \ \leq \ \eta.
\end{equation}
There are absolute constants $C,c,C_{\mathrm{rec}}>0$ such that if
\begin{equation}\label{eq:abstractsamples}
q \ \geq \ C K_L^2\frac{r^2}{\varepsilon^2}\log^2\left(2+K_L^2\frac{r^2}{\varepsilon^2}\right)\log(2D_L),
\end{equation}
then, with probability at least $1-D_L^{-c}$,
\begin{equation}\label{eq:abstractnoisyrecovery}
\left\|f^*-f\right\|_{L^2(X,\mu)} \ \leq \ C_{\mathrm{rec}}\left(\varepsilon\left\|f\right\|_{L^2(X,\mu)}+\eta\right).
\end{equation}
The same event works for all observations whose errors satisfy the prescribed bound. In particular, exact observations allow $\eta=0$.
\end{theorem}

\begin{proof}
Write
\begin{equation}
f \ = \ \sum_{\lambda_j\leq L}a_je_j.
\end{equation}
Parseval's identity and the Fourier Ratio hypothesis give
\begin{equation}
\left\|a\right\|_2 \ = \ \left\|f\right\|_{L^2}, \qquad \left\|a\right\|_1 \ \leq \ r\left\|a\right\|_2.
\end{equation}
Set
\begin{equation}
u \ = \ \frac{r^2}{\varepsilon^2}, \qquad S \ = \ \min\left\{D_L,\left\lceil u\right\rceil\right\}, \qquad t \ = \ \min\{D_L,2S\}.
\end{equation}
Here $r\geq1$, since the $\ell^1$ norm of a nonzero vector is at least its $\ell^2$ norm. Thus $u>4$ and $t\leq2\lceil u\rceil\leq(5/2)u$.

For a coefficient vector $a$, let $\sigma_S(a)_1$ be the smallest $\ell^1$ error among all approximations supported on at most $S$ coordinates. If $S<D_L$, then
\begin{equation}\label{eq:FRbesttermapproximation}
\frac{\sigma_S(a)_1}{\sqrt S} \ \leq \ \frac{\left\|a\right\|_1}{\sqrt S} \ \leq \ \varepsilon\left\|a\right\|_2.
\end{equation}
If $S=D_L$, the approximation error is zero.

Form the sampling matrix
\begin{equation}
E_{ij} \ = \ \frac{1}{\sqrt q}e_j(x_i), \qquad \lambda_j\leq L.
\end{equation}
Its rows are independent, its columns are orthonormal in expectation, and its entries are bounded in magnitude by $K_L/\sqrt q$. Fix a sufficiently small absolute restricted-isometry constant, for example $1/4$. The bounded-orthonormal-system estimate in \cite[Theorem~2.3]{BrugiapagliaDirksenJungRauhut} gives the restricted-isometry property at order $t$ when
\begin{equation}\label{eq:BOSsamplecondition}
q \ \geq \ C K_L^2t\log^2\left(2+K_L^2t\right)\log(2D_L).
\end{equation}
Its exponential probability bound implies probability at least $1-D_L^{-c}$ after adjusting the absolute constants. Since $t\leq(5/2)u$, condition~\eqref{eq:abstractsamples} implies~\eqref{eq:BOSsamplecondition}.

The coefficient form of the optimization problem has data $q^{-1/2}(g_i)_{i=1}^q$ and residual radius $\eta$. If $2S<D_L$, the restricted-isometry property at order $2S$ gives the stable basis-pursuit estimate
\begin{equation}
\left\|a^*-a\right\|_2 \ \leq \ C_0\frac{\sigma_S(a)_1}{\sqrt S}+C_1\eta.
\end{equation}
Combining this with~\eqref{eq:FRbesttermapproximation} proves the required bound in this case.

If $2S\geq D_L$, then $t=D_L$ and the restricted-isometry property controls every coefficient vector. Both $a$ and $a^*$ are feasible, so
\begin{equation}
\left\|E(a^*-a)\right\|_2 \ \leq \ 2\eta.
\end{equation}
The lower restricted-isometry inequality therefore gives
\begin{equation}
\left\|a^*-a\right\|_2 \ \leq \ \frac{2\eta}{\sqrt{1-1/4}}.
\end{equation}
This also implies~\eqref{eq:abstractnoisyrecovery}. Finally, Parseval's identity converts coefficient error into function error.
\end{proof}

The radius $\eta$ is part of the input to the reconstruction problem. In an application it must come from a known noise bound or from a known bound on the energy of the signal. Exact observations require no such additional information because one may impose exact agreement with the samples.

\begin{corollary}\label{cor:spectralregularityrecovery}
Assume the hypotheses of Theorem~\ref{thm:spectralrecovery}, and fix $s>0$. For nonzero $f\in V_L$, define
\begin{equation}
r_{L,s}(f) \ = \ \frac{\left|\int_Xf\,d\mu\right|}{\left\|f\right\|_{L^2}}+\left(\sum_{0<\lambda_j\leq L}\lambda_j^{-s}\right)^{1/2}\frac{\left\|A^{s/2}f\right\|_{L^2}}{\left\|f\right\|_{L^2}}.
\end{equation}
For observations and reconstruction radius as in that theorem, the estimate
\begin{equation}
\left\|f^*-f\right\|_{L^2(X,\mu)} \ \leq \ C_{\mathrm{rec}}\left(\varepsilon\left\|f\right\|_{L^2(X,\mu)}+\eta\right)
\end{equation}
holds with probability at least $1-D_L^{-c}$ whenever
\begin{equation}
q \ \geq \ C K_L^2\frac{r_{L,s}(f)^2}{\varepsilon^2}\log^2\left(2+K_L^2\frac{r_{L,s}(f)^2}{\varepsilon^2}\right)\log(2D_L).
\end{equation}
\end{corollary}

\begin{proof}
Theorem~\ref{thm:spectral} gives $\operatorname{FR}_L(f)\leq r_{L,s}(f)$. Apply Theorem~\ref{thm:spectralrecovery}.
\end{proof}

% ========================= Recovery on Compact Resistance Spaces =========================
\subsection{Recovery on Compact Resistance Spaces}
On a compact measured resistance space, the resistance inequality controls the oscillation of each eigenfunction. The eigenvalue equation then expresses its energy in terms of its eigenvalue.

\begin{proposition}\label{prop:resistancecoherence}
Let $(X,R,\mu,\mathcal E,\mathcal F)$ be a compact measured resistance space with resistance diameter $D_R$, and let $A=-\Delta_\mu$ have the eigenbasis of Theorem~\ref{thm:reslap}. For every nonconstant eigenfunction,
\begin{equation}
\left\|e_j\right\|_{L^\infty(X)} \ \leq \ \sqrt{D_R\lambda_j}.
\end{equation}
Consequently,
\begin{equation}
K_L^2 \ \leq \ \max\{1,D_RL\}.
\end{equation}
\end{proposition}

\begin{proof}
For $j\geq1$, orthogonality to the constant eigenfunction gives $\int_Xe_j\,d\mu=0$. Hence, for every $x\in X$,
\begin{equation}
e_j(x) \ = \ \int_X\left(e_j(x)-e_j(y)\right)\,d\mu(y).
\end{equation}
The resistance inequality gives
\begin{equation}
\left|e_j(x)-e_j(y)\right| \ \leq \ R(x,y)^{1/2}\mathcal E(e_j,e_j)^{1/2} \ = \ R(x,y)^{1/2}\lambda_j^{1/2}.
\end{equation}
Since $R(x,y)\leq D_R$ and $\mu(X)=1$, integration yields $\left|e_j(x)\right|\leq\sqrt{D_R\lambda_j}$. The constant eigenfunction has supremum norm $1$.
\end{proof}

\begin{theorem}\label{thm:resistancegeometricrecovery}
Let $(X,R,\mu,\mathcal E,\mathcal F)$ be a compact measured resistance space with resistance diameter $D_R$ and resistance Laplacian $A=-\Delta_\mu$. Define
\begin{equation}
C_R \ = \ \left(\frac12\int_X\int_XR(x,y)\,d\mu(x)\,d\mu(y)\right)^{1/2}, \qquad M_L \ = \ \max\{1,D_RL\},
\end{equation}
and, for nonzero $f\in V_L$, define
\begin{equation}
r_R(f) \ = \ \frac{\left|\int_Xf\,d\mu\right|}{\left\|f\right\|_{L^2}}+C_R\frac{\mathcal E(f,f)^{1/2}}{\left\|f\right\|_{L^2}}.
\end{equation}
Let $0<\varepsilon<1/2$, and use the observations and prescribed radius $\eta$ from Theorem~\ref{thm:spectralrecovery}. If
\begin{equation}\label{eq:resistancesamples}
q \ \geq \ C M_L\frac{r_R(f)^2}{\varepsilon^2}\log^2\left(2+M_L\frac{r_R(f)^2}{\varepsilon^2}\right)\log(2D_L),
\end{equation}
then, with probability at least $1-D_L^{-c}$,
\begin{equation}
\left\|f^*-f\right\|_{L^2(X,\mu)} \ \leq \ C_{\mathrm{rec}}\left(\varepsilon\left\|f\right\|_{L^2(X,\mu)}+\eta\right).
\end{equation}
Moreover,
\begin{equation}
C_R \ \leq \ \left(\frac{D_R}{2}\right)^{1/2}, \qquad r_R(f) \ \leq \ 1+\left(\frac{D_R}{2}\right)^{1/2}\frac{\mathcal E(f,f)^{1/2}}{\left\|f\right\|_{L^2}}.
\end{equation}
The right side of the second inequality may replace $r_R(f)$ in the sample-size condition. Exact observations allow $\eta=0$.
\end{theorem}

\begin{proof}
Corollary~\ref{cor:uniformresistance} gives $\operatorname{FR}_L(f)\leq r_R(f)$, while Proposition~\ref{prop:resistancecoherence} gives $K_L^2\leq M_L$. Substitute both bounds into Theorem~\ref{thm:spectralrecovery}. The remaining estimates follow from $R(x,y)\leq D_R$, $\mu(X)=1$, and Cauchy--Schwarz.
\end{proof}

The Fourier Ratio constant $C_R$ is independent of the cutoff. The coherence bound used in~\eqref{eq:resistancesamples} still grows with the cutoff, however, and its size has an important consequence for the interpretation of that theorem.

\begin{proposition}\label{prop:diametersamplinglimitation}
For every spectral cutoff $L>0$ on a compact measured resistance space,
\begin{equation}
D_L-1 \ \leq \ L\sum_{j\geq1}\frac1{\lambda_j} \ = \ LC_R^2 \ \leq \ \frac{D_RL}{2},
\end{equation}
where the sum is over the positive eigenvalues and has finitely many terms when $X$ is finite. If $D_L\geq2$, then
\begin{equation}
\max\{1,D_RL\} \ \geq \ 2(D_L-1) \ \geq \ D_L.
\end{equation}
\end{proposition}

\begin{proof}
Each positive eigenvalue below $L$ satisfies $1\leq L/\lambda_j$. Sum this inequality and apply Theorem~\ref{thm:resistancetrace} and $C_R^2\leq D_R/2$.
\end{proof}

Thus the purely geometric budget in~\eqref{eq:resistancesamples} already contains a factor at least as large as the dimension $D_L$ of the spectral subspace before the accuracy and logarithmic factors are included. It provides a sufficient oversampling bound. To obtain a compressed-sampling conclusion, we must retain the actual $K_L$ in Theorem~\ref{thm:spectralrecovery} and estimate it more accurately. The following example gives a case in which this can be done explicitly.

\subsection{An Interval Example}
The counting-function illustration in Figure~\ref{fig:weyl-counting} uses the Neumann Laplacian on $[0,\pi]$, with positive eigenvalues $j^2$. The interval $[0,1]$ used below differs only by the rescaling $y=\pi x$, so its positive Neumann eigenvalues are $\pi^2j^2$. With probability-normalized Lebesgue measures, the corresponding orthonormal eigenfunctions have the same $L^\infty$ bound.

\begin{proposition}\label{prop:intervalsampling}
Let $X=[0,1]$, let $\mu$ be Lebesgue measure, and take
\begin{equation}
\mathcal F \ = \ H^1([0,1]), \qquad \mathcal E(f,f) \ = \ \int_0^1\left|f'(x)\right|^2\,dx.
\end{equation}
The associated operator is the Neumann Laplacian. Its resistance metric and resistance constant are
\begin{equation}
R(x,y) \ = \ \left|x-y\right|, \qquad C_R^2 \ = \ \frac16.
\end{equation}
For the orthonormal eigenbasis
\begin{equation}
e_0(x)=1, \qquad e_j(x)=\sqrt2\cos(\pi jx), \qquad \lambda_j=\pi^2j^2 \quad (j\geq1),
\end{equation}
one has $K_L^2=2$ whenever $L\geq\pi^2$, and
\begin{equation}
D_L \ = \ 1+\left\lfloor\frac{\sqrt L}{\pi}\right\rfloor.
\end{equation}
Suppose that a nonzero $f\in V_L$ satisfies
\begin{equation}
\left\|f'\right\|_{L^2([0,1])} \ \leq \ H\left\|f\right\|_{L^2([0,1])},
\end{equation}
and set $r_H=1+H/\sqrt6$. For exact observations, stable reconstruction with relative error at most $C_{\mathrm{rec}}\varepsilon$ holds with probability at least $1-D_L^{-c}$ if
\begin{equation}
q \ \geq \ C\frac{r_H^2}{\varepsilon^2}\log^2\left(2+\frac{r_H^2}{\varepsilon^2}\right)\log(2D_L).
\end{equation}
For fixed $H$ and $\varepsilon$, this budget grows logarithmically in $D_L$.
\end{proposition}

\begin{proof}
For $x<y$, the fundamental theorem of calculus and Cauchy--Schwarz give
\begin{equation}
\left|f(y)-f(x)\right|^2 \ \leq \ (y-x)\int_x^y\left|f'(t)\right|^2\,dt \ \leq \ (y-x)\mathcal E(f,f).
\end{equation}
Equality in the defining resistance supremum is attained by a function that is constant on $[0,x]$, linear on $[x,y]$, and constant on $[y,1]$. Thus $R(x,y)=y-x$. Consequently,
\begin{equation}
C_R^2 \ = \ \frac12\int_0^1\int_0^1\left|x-y\right|\,dx\,dy \ = \ \frac16.
\end{equation}
The eigenfunctions and eigenvalues follow from the Neumann boundary conditions. Their supremum norms give $K_L^2=2$. Corollary~\ref{cor:uniformresistance} gives $\operatorname{FR}_L(f)\leq r_H$, and Theorem~\ref{thm:spectralrecovery} with $\eta=0$ proves the sampling assertion after absorbing the fixed factors into $C$.
\end{proof}

% ========================= Sensor Stability =========================
\subsection{Sensor Stability in the Resistance Metric}
The resistance inequality also controls the effect of moving a sensor. This part of the argument is deterministic and requires no assumption on the distribution of the sensor locations.

\begin{theorem}\label{thm:resistancesensor}
Let $f\in\mathcal F$, and suppose that intended sensor locations $x_1,\ldots,x_q$ and actual locations $y_1,\ldots,y_q$ satisfy
\begin{equation}
R(x_i,y_i) \ \leq \ \delta
\end{equation}
for every $i$. Then
\begin{equation}\label{eq:sensorerror}
\left(\frac1q\sum_{i=1}^q\left|f(y_i)-f(x_i)\right|^2\right)^{1/2} \ \leq \ \delta^{1/2}\mathcal E(f,f)^{1/2}.
\end{equation}
If $\mathcal E(f,f)>0$, the placement error is at most $\varepsilon\left\|f\right\|_{L^2(X,\mu)}$ whenever
\begin{equation}
\delta \ \leq \ \varepsilon^2\frac{\left\|f\right\|_{L^2(X,\mu)}^2}{\mathcal E(f,f)}.
\end{equation}
If $\mathcal E(f,f)=0$, then $f$ is constant and the placement error is zero. For $f\in V_L$ and $L>0$, the sufficient condition
\begin{equation}
\delta \ \leq \ \frac{\varepsilon^2}{L}
\end{equation}
guarantees the same relative error bound.
\end{theorem}

\begin{proof}
For each $i$, the resistance inequality gives
\begin{equation}
\left|f(y_i)-f(x_i)\right|^2 \ \leq \ R(x_i,y_i)\mathcal E(f,f) \ \leq \ \delta\mathcal E(f,f).
\end{equation}
Average over the measurements and take square roots. For $f\in V_L$, the spectral expansion gives
\begin{equation}
\mathcal E(f,f) \ = \ \sum_{\lambda_j\leq L}\lambda_j\left|\widehat{f}(j)\right|^2 \ \leq \ L\left\|f\right\|_{L^2}^2,
\end{equation}
which proves the last assertion.
\end{proof}

\begin{theorem}\label{thm:perturbedresistancerecovery}
Let $0<\varepsilon<1/2$, and let $f\in V_L$ be nonzero with $\operatorname{FR}_L(f)\leq r$. Choose intended locations $x_1,\ldots,x_q$ independently according to $\mu$, and suppose that the actual locations satisfy $R(x_i,y_i)\leq\delta$. Set
\begin{equation}
g_i \ = \ f(y_i), \qquad \eta \ = \ \delta^{1/2}\mathcal E(f,f)^{1/2}.
\end{equation}
Use the observations $g_i$ and radius $\eta$ in the reconstruction problem of Theorem~\ref{thm:spectralrecovery}. If
\begin{equation}
q \ \geq \ C K_L^2\frac{r^2}{\varepsilon^2}\log^2\left(2+K_L^2\frac{r^2}{\varepsilon^2}\right)\log(2D_L),
\end{equation}
then, with probability at least $1-D_L^{-c}$,
\begin{equation}\label{eq:perturbedrecoveryerror}
\left\|f^*-f\right\|_{L^2(X,\mu)} \ \leq \ C_{\mathrm{rec}}\left(\varepsilon\left\|f\right\|_{L^2(X,\mu)}+\delta^{1/2}\mathcal E(f,f)^{1/2}\right).
\end{equation}
In particular, if $L>0$ and $\delta\leq\varepsilon^2/L$, then
\begin{equation}
\left\|f^*-f\right\|_{L^2(X,\mu)} \ \leq \ C'_{\mathrm{rec}}\varepsilon\left\|f\right\|_{L^2(X,\mu)}.
\end{equation}
A known upper bound for $\mathcal E(f,f)^{1/2}$ may be used in the reconstruction radius, with the same replacement in the error estimate.
\end{theorem}

\begin{proof}
Write $\xi_i=f(y_i)-f(x_i)$. Theorem~\ref{thm:resistancesensor} gives
\begin{equation}
\left(\frac1q\sum_{i=1}^q\left|\xi_i\right|^2\right)^{1/2} \ \leq \ \eta.
\end{equation}
Theorem~\ref{thm:spectralrecovery} therefore applies with this prescribed radius and gives~\eqref{eq:perturbedrecoveryerror}. Its sampling event depends on the intended locations and works for every perturbation satisfying the resistance bound. The final assertion follows from $\mathcal E(f,f)\leq L\left\|f\right\|_{L^2}^2$.
\end{proof}

\begin{corollary}\label{cor:resistancecells}
Let $Q_1,\ldots,Q_K\subset X$ satisfy
\begin{equation}
\operatorname{diam}_R(Q_k) \ \leq \ h
\end{equation}
for every $k$. If each intended sensor and its actual location lie in the same cell, then
\begin{equation}
\left|f(y_i)-f(x_i)\right| \ \leq \ h^{1/2}\mathcal E(f,f)^{1/2}.
\end{equation}
The conclusions of Theorems~\ref{thm:resistancesensor} and~\ref{thm:perturbedresistancerecovery} hold with $\delta=h$.
\end{corollary}

\begin{proof}
If $x_i,y_i\in Q_k$, then $R(x_i,y_i)\leq\operatorname{diam}_R(Q_k)\leq h$. Apply the resistance inequality and the two theorems.
\end{proof}

% ========================= Fractals =========================
\section{Fractals}
We now specialize the preceding results to self-similar fractals. There are two ways to estimate the reciprocal eigenvalue sum in this setting. A modified Weyl law describes its convergence through the spectral dimension. Kasue's resistance identity, recalled in Theorem~\ref{thm:resistancetrace}, gives its value directly from the average resistance distance \cite{Kasue}. Our application of this identity to the Fourier Ratio does not require a counting law.

% ========================= Spectral Dimension =========================
\subsection{Spectral Dimension}
For many self-similar fractals the eigenvalue counting function has the form of a power $\Lambda^{d_s/2}$ multiplied by a bounded log-periodic factor \cite{Riane,KigamiLapidus,KumagaiSpectralDimension}. The exponent $d_s$ is called the spectral dimension. Under the lower bound on the oscillatory factor stated below, the sublinear counting result forces $d_s<2$.

\begin{theorem}\label{thm:lowdimension}
Let $(X,R,\mu,\mathcal E,\mathcal F)$ be a compact measured resistance space with $\mu(X) \ = \ 1$ and resistance Laplacian $A \ = \ -\Delta_\mu$. Suppose that for some $d_s>0$ the eigenvalue counting function satisfies a modified Weyl law
\begin{equation}
N_A(\Lambda) \ = \ \left(G(\log\Lambda)+o(1)\right)\Lambda^{d_s/2}
\end{equation}
as $\Lambda\rightarrow\infty$, where $G$ is bounded above and bounded away from zero. Then
\begin{equation}d_s \ < \ 2.\end{equation}
Consequently Theorem~\ref{thm:counting} with $s \ = \ 1$ gives a cutoff-independent bound for the reciprocal eigenvalue sum.
\end{theorem}

\begin{proof}
Corollary~\ref{cor:sublinearcounting} tells us that
\begin{equation}\frac{N_A(\Lambda)}{\Lambda} \ \longrightarrow \ 0.\end{equation}
If $d_s>2$, the modified Weyl law and the positive lower bound on $G$ make this quotient tend to infinity. If $d_s \ = \ 2$, then the quotient is $G(\log\Lambda)+o(1)$ and again cannot tend to zero. Therefore $d_s<2$. The counting estimate has exponent $d_s/2<1$, so Theorem~\ref{thm:counting} with $s \ = \ 1$ also recovers convergence of the reciprocal eigenvalue sum from the spectral side.
\end{proof}

We do not need to assume $d_s<2$ to obtain the Fourier Ratio bound on a compact resistance space. The resistance identity already supplies that bound. The theorem shows instead that a modified Weyl law with its leading factor bounded away from zero must have $d_s<2$.

% ========================= The Sierpinski Gasket =========================
\subsection{The Sierpinski Gasket}
We use the standard Sierpinski gasket as a concrete example. Its spectral decimation describes the counting exponent, while its resistance form supplies the geometric constant in the Fourier Ratio estimate. Throughout this section we use the operator associated with the full energy domain, which is the Neumann Laplacian.

\begin{figure}[H]
\centering
\def\trianglewidth{2cm}
\pgfdeclarelindenmayersystem{Sierpinski triangle}{
    \symbol{X}{\pgflsystemdrawforward}
    \symbol{Y}{\pgflsystemdrawforward}
    \rule{X -> X-Y+X+Y-X}
    \rule{Y -> YY}
}
\foreach \level in {0,...,6}{
\tikzset{
    l-system={step=\trianglewidth/(2^\level), order=\level, angle=-120}
}
\begin{tikzpicture}
    \fill [black] (0,0) -- ++(0:\trianglewidth) -- ++(120:\trianglewidth) -- cycle;
    \draw [draw=none] (0,0) l-system
    [l-system={Sierpinski triangle, axiom=X},fill=white];
\end{tikzpicture}
}
\caption{Finite approximations of the Sierpiński gasket from levels $0$ through $6$, illustrating the recursive self-similar construction obtained by successively removing central triangular regions.
}
\label{fig:sierpinski-approximations}
\end{figure}
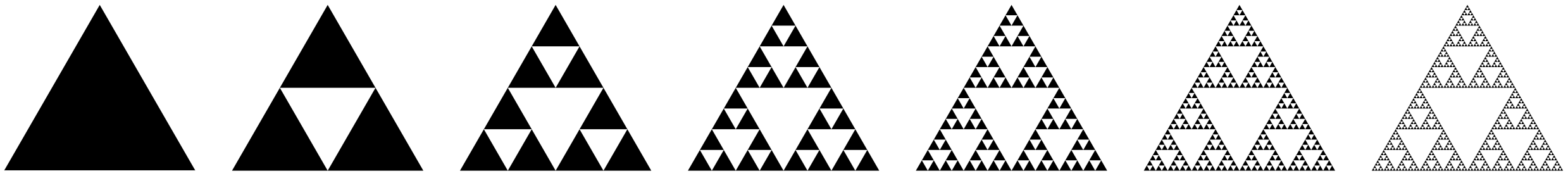

\begin{proposition}\label{prop:sgdimension}
For the standard Neumann Laplacian on the Sierpinski gasket,
\begin{equation}N_A(\Lambda) \ \asymp \ \Lambda^{\log 3/\log 5}.\end{equation}
The more precise counting law has a bounded log-periodic factor \cite{FukushimaShima,StrichartzSpectral}. In particular,
\begin{equation}\frac{d_s}{2} \ = \ \frac{\log 3}{\log 5}<1, \qquad d_s \ = \ \frac{2\log 3}{\log 5}<2.\end{equation}
\end{proposition}

\begin{proof}
Spectral decimation gives an eigenvalue scale $5^m$ and a counting scale $3^m$, including multiplicities \cite{FukushimaShima}. Thus
\begin{equation}
N_A(5^m) \ \asymp \ 3^m \ = \ \left(5^m\right)^{\log 3/\log 5}.
\end{equation}
If $5^m\leq\Lambda<5^{m+1}$, monotonicity of $N_A$ gives
\begin{equation}
N_A(5^m) \ \leq \ N_A(\Lambda) \ \leq \ N_A(5^{m+1}).
\end{equation}
The first counting estimate follows, and comparison with $\Lambda^{d_s/2}$ identifies $d_s$. The log-periodic refinement requires the more precise spectral asymptotics in the cited works and does not follow from this comparison alone.
\end{proof}

\begin{corollary}\label{cor:sg}
Let $SG$ be the standard Sierpinski gasket with its standard resistance form and self-similar probability measure. For every spectral cutoff $L$ and every nonzero $f\in V_L$,
\begin{equation}
\operatorname{FR}_L(f) \ \leq \ \frac{\left|\int_{SG}f\,d\mu\right|}{\left\|f\right\|_{L^2(SG,\mu)}}+C_{SG}\frac{\mathcal E(f,f)^{1/2}}{\left\|f\right\|_{L^2(SG,\mu)}} \ \leq \ 1+C_{SG}\frac{\mathcal E(f,f)^{1/2}}{\left\|f\right\|_{L^2(SG,\mu)}},
\end{equation}
where
\begin{equation}
C_{SG} \ = \ \left(\frac{1}{2}\int_{SG}\int_{SG}R(x,y)\,d\mu(x)\,d\mu(y)\right)^{1/2} \ \leq \ \left(\frac{D_R}{2}\right)^{1/2}.
\end{equation}
In particular, $C_{SG}$ does not depend on $L$.
\end{corollary}

\begin{proof}
Apply Corollary~\ref{cor:uniformresistance} to the standard resistance form on $SG$.
\end{proof}

The estimate controls the Fourier Ratio across all spectral cutoffs for families with a fixed upper bound on
\begin{equation}
\frac{\mathcal E(f,f)}{\left\|f\right\|_{L^2(SG,\mu)}^2}.
\end{equation}
The energy term is essential, since it may increase as more of the spectrum is admitted.

% ========================= Further Fractal Examples =========================
\subsection{Further Fractal Examples}
The Sierpinski tetrahedron and the Minkowski curve give two further examples \cite{Riane,RianeDavidTetrahedron,RianeDavidMinkowski}.

\begin{proposition}\label{prop:otherfractals}
For the standard Neumann Laplacians on the Sierpinski tetrahedron $ST$ and the Minkowski curve $MC$,
\begin{equation}
N_{ST}(\Lambda) \ \asymp \ \Lambda^{\log 4/\log 6}, \qquad N_{MC}(\Lambda) \ \asymp \ \Lambda^{1/2}.
\end{equation}
Consequently,
\begin{equation}
d_s(ST) \ = \ \frac{2\log 4}{\log 6}<2, \qquad d_s(MC) \ = \ 1<2.
\end{equation}
\end{proposition}

\begin{proof}
For the Sierpinski tetrahedron spectral decimation gives an eigenvalue scale $6^m$ and a counting scale $4^m$ \cite{RianeDavidTetrahedron}. The comparison used for the gasket therefore gives the exponent $\log 4/\log 6$.

For the Minkowski curve the corresponding scales are $64^m$ and $8^m$ \cite{RianeDavidMinkowski}. Thus the counting exponent is
\begin{equation}
\frac{\log 8}{\log 64} \ = \ \frac12.
\end{equation}
Comparing these exponents with $d_s/2$ gives the stated spectral dimensions. These scale comparisons establish the two-sided counting bounds. More precise oscillatory asymptotics require the spectral results in the cited works.
\end{proof}

\begin{corollary}\label{cor:otherfractals}
Let $X$ be either the Sierpinski tetrahedron or the Minkowski curve, equipped with its standard resistance form and self-similar probability measure. For every spectral cutoff $L$ and every nonzero $f\in V_L$,
\begin{equation}
\operatorname{FR}_L(f) \ \leq \ 1+C_X\frac{\mathcal E(f,f)^{1/2}}{\left\|f\right\|_{L^2(X,\mu)}},
\end{equation}
where
\begin{equation}
C_X \ = \ \left(\frac12\int_X\int_XR(x,y)\,d\mu(x)\,d\mu(y)\right)^{1/2}
\end{equation}
is independent of $L$.
\end{corollary}

\begin{proof}
Apply Corollary~\ref{cor:uniformresistance} to the corresponding compact measured resistance space.
\end{proof}

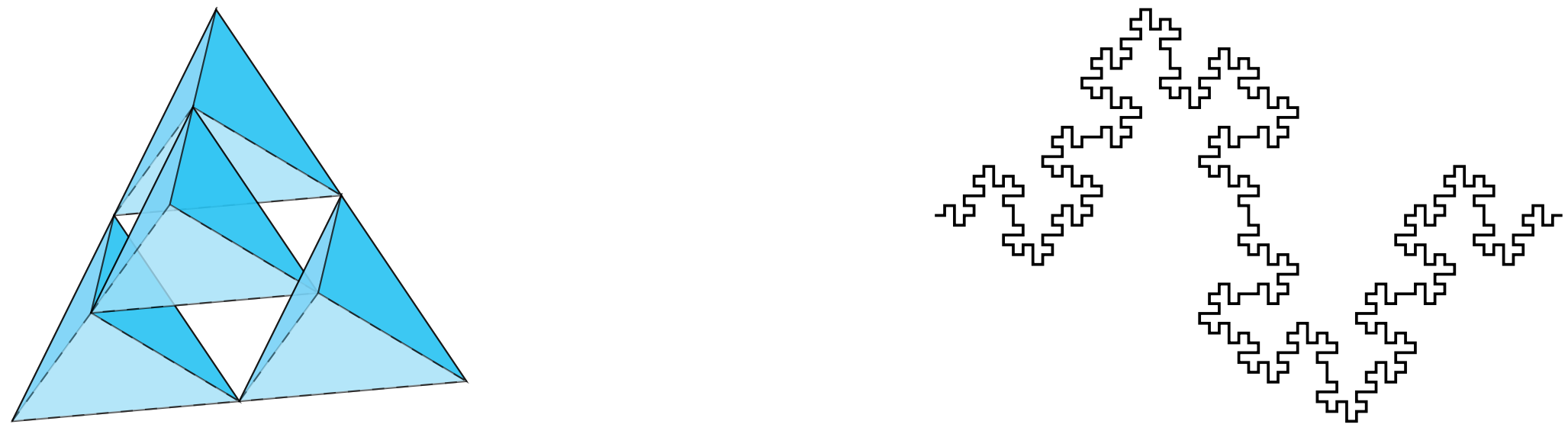
\begin{figure}[H]
\centering
\begin{subfigure}[c]{0.47\textwidth}
\centering

\tdplotsetmaincoords{60}{110}

\begin{tikzpicture}[tdplot_main_coords, scale=2.0]

    \coordinate (A) at (0, 0, 1.633);
    \coordinate (B) at (1, -0.577, 0);
    \coordinate (C) at (-1, -0.577, 0);
    \coordinate (D) at (0, 1.155, 0);

    \newcommand{\drawtet}[2]{%
        \begin{scope}[shift={#1}, scale=#2]
            \draw[fill=cyan!30, opacity=0.7, thin]
                (0,0,1.633) -- (1,-0.577,0) -- (-1,-0.577,0) -- cycle;

            \draw[fill=cyan!50, opacity=0.7, thin]
                (0,0,1.633) -- (1,-0.577,0) -- (0,1.155,0) -- cycle;

            \draw[fill=cyan!70, opacity=0.7, thin]
                (0,0,1.633) -- (-1,-0.577,0) -- (0,1.155,0) -- cycle;

            \draw[fill=cyan!10, opacity=0.4, dashed]
                (1,-0.577,0) -- (-1,-0.577,0) -- (0,1.155,0) -- cycle;
        \end{scope}
    }

    \drawtet{(0,0,0.8165)}{0.5}
    \drawtet{(0.5,-0.2885,0)}{0.5}
    \drawtet{(-0.5,-0.2885,0)}{0.5}
    \drawtet{(0,0.5775,0)}{0.5}

\end{tikzpicture}

\caption{Sierpiński tetrahedron}
\label{fig:sierpinski-tetrahedron}

\end{subfigure}
\hfill
\begin{subfigure}[c]{0.47\textwidth}
\centering

\pgfdeclarelindenmayersystem{Minkowski}{
    \rule{F -> F+F-F-FF+F+F-F}
}

\begin{tikzpicture}

    \draw[
        line width=0.7pt,
        l-system={
            Minkowski,
            axiom=F,
            order=3,
            angle=90,
            step=0.085cm
        }
    ]
    l-system;

\end{tikzpicture}

\caption{Minkowski curve}
\label{fig:minkowski-curve}

\end{subfigure}

\caption{Finite fractal geometries used as examples of self-similar resistive networks: a Sierpiński tetrahedron and a finite approximation of the Minkowski curve.}
\label{fig:fractal-networks}

\end{figure}

% ========================= Reconstruction on Fractals =========================
\subsection{Reconstruction on Fractals}
The Fourier Ratio estimates now enter the sampling theorem. The resistance constant controls the coefficient complexity. The additional quantity $K_L$ controls the size of the eigenfunctions at the sampling points. Keeping track of both quantities is necessary to determine whether the sufficient number of samples is smaller than the dimension of $V_L$.

\begin{corollary}\label{cor:sgrecovery}
Let $SG$ be the standard Sierpinski gasket with its standard resistance form and self-similar probability measure, and let $C_{SG}$ be the constant from Corollary~\ref{cor:sg}. Fix $0<\varepsilon<1/2$. For nonzero $f\in V_L$, define
\begin{equation}
r_{SG}(f) \ = \ \frac{\left|\int_{SG}f\,d\mu\right|}{\left\|f\right\|_{L^2(SG,\mu)}}+C_{SG}\frac{\mathcal E(f,f)^{1/2}}{\left\|f\right\|_{L^2(SG,\mu)}},
\qquad B_L \ = \ \max\left\{1,D_RL\right\}.
\end{equation}
Choose $x_1,\ldots,x_q$ independently according to $\mu$ and suppose that
\begin{equation}
q \ \geq \ C B_L\frac{r_{SG}(f)^2}{\varepsilon^2}\log^2\left(2+\frac{B_Lr_{SG}(f)^2}{\varepsilon^2}\right)\log(2D_L).
\end{equation}
For sample data with noise bounded by $\eta$ in the normalized sampled $\ell^2$ norm, let $f^*$ be the reconstruction from Theorem~\ref{thm:spectralrecovery} with constraint radius $\eta$. Then with probability at least $1-D_L^{-c}$,
\begin{equation}
\left\|f^*-f\right\|_{L^2(SG,\mu)} \ \leq \ C_{\mathrm{rec}}\left(\varepsilon\left\|f\right\|_{L^2(SG,\mu)}+\eta\right).
\end{equation}
For exact sample data one may take $\eta=0$.
\end{corollary}

\begin{proof}
Corollary~\ref{cor:sg} gives $\operatorname{FR}_L(f)\leq r_{SG}(f)$, while Proposition~\ref{prop:resistancecoherence} gives $K_L^2\leq B_L$. Substitute these estimates into Theorem~\ref{thm:spectralrecovery}.
\end{proof}

\begin{corollary}\label{cor:otherfractalrecovery}
Let $X$ be either the Sierpinski tetrahedron or the Minkowski curve with the resistance form and probability measure used above, and let $C_X$ be the constant from Corollary~\ref{cor:otherfractals}. Fix $0<\varepsilon<1/2$. For nonzero $f\in V_L$, define
\begin{equation}
r_X(f) \ = \ \frac{\left|\int_Xf\,d\mu\right|}{\left\|f\right\|_{L^2(X,\mu)}}+C_X\frac{\mathcal E(f,f)^{1/2}}{\left\|f\right\|_{L^2(X,\mu)}},
\qquad B_L \ = \ \max\left\{1,D_RL\right\}.
\end{equation}
Choose $x_1,\ldots,x_q$ independently according to $\mu$ and suppose that
\begin{equation}
q \ \geq \ C B_L\frac{r_X(f)^2}{\varepsilon^2}\log^2\left(2+\frac{B_Lr_X(f)^2}{\varepsilon^2}\right)\log(2D_L).
\end{equation}
For sample data with noise bounded by $\eta$ in the normalized sampled $\ell^2$ norm, let $f^*$ be the reconstruction from Theorem~\ref{thm:spectralrecovery} with constraint radius $\eta$. Then with probability at least $1-D_L^{-c}$,
\begin{equation}
\left\|f^*-f\right\|_{L^2(X,\mu)} \ \leq \ C_{\mathrm{rec}}\left(\varepsilon\left\|f\right\|_{L^2(X,\mu)}+\eta\right).
\end{equation}
For exact sample data one may take $\eta=0$.
\end{corollary}

\begin{proof}
Corollary~\ref{cor:uniformresistance} gives $\operatorname{FR}_L(f)\leq r_X(f)$. Use Proposition~\ref{prop:resistancecoherence} and Theorem~\ref{thm:spectralrecovery} as above.
\end{proof}

These corollaries give geometric sufficient conditions for stable reconstruction, but the diameter estimate alone does not establish a saving in the number of samples. On each example above, $D_L$ grows like $L^{d_s/2}$ with $d_s/2<1$, whereas $B_L$ grows linearly in $L$. The sufficient sample count obtained by substituting $B_L$ therefore eventually exceeds $D_L$ for fixed $\varepsilon$. A smaller sample count requires a sharper estimate for the actual $K_L$, together with control of the energy quotient in $r_X(f)$. Theorem~\ref{thm:spectralrecovery} permits these sharper estimates to be used directly.

Finally, Theorems~\ref{thm:resistancesensor} and~\ref{thm:perturbedresistancerecovery} apply to all three examples. They control changes in measured values by resistance distance and energy, and hence give the corresponding reconstruction estimates when the sampling locations are perturbed.

\section{Electrical Networks}
Finite resistive electrical networks give a direct physical application of the resistance-space theory. The resistance metric is effective electrical resistance, and the energy of a voltage potential is the power dissipated in the network \cite{DoyleSnell,KleinRandic,LyonsPeres}. The associated operator is a weighted graph Laplacian. Thus the general estimates above have an electrical interpretation in terms of voltage, power, and effective resistance.

Sampling and reconstruction of functions on graphs have an established theory. Chen, Varma, Sandryhaila, and Kova\v{c}evi\'c develop sampling theory for graph signals \cite{ChenVarmaSandryhailaKovacevic}, and Puy, Tremblay, Gribonval, and Vandergheynst study random sampling and stable reconstruction of bandlimited graph signals \cite{PuyTremblayGribonvalVandergheynst}. Our purpose here is to express the Fourier Ratio estimates in electrical quantities and to account for the additional constraint imposed by a current source. The current-voltage equation divides each nonconstant spectral coefficient of the current by its Laplacian eigenvalue. We use this elementary observation to control the part of the voltage above a spectral cutoff.

The geometric sampling bounds remain sufficient bounds. As explained after Theorem~\ref{thm:resistancegeometricrecovery}, the general diameter estimate does not by itself guarantee reconstruction from fewer samples than the dimension of the spectral subspace. To obtain such a conclusion for a particular network, one must retain sharper information about its eigenfunctions and the Fourier Ratio of the voltage.

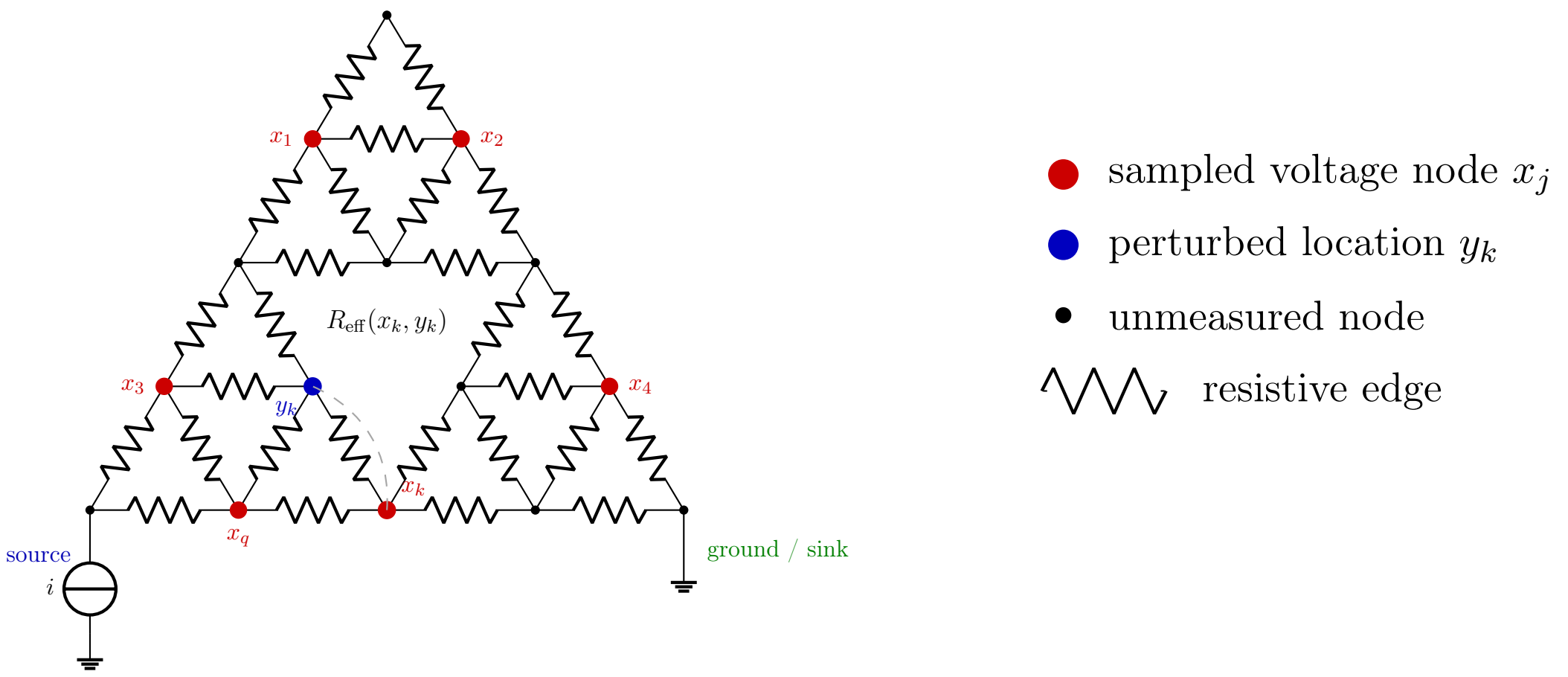
\begin{figure}[H]
\centering
\begin{circuitikz}[american resistors, scale=0.55, transform shape]
% Sierpinski-gasket vertices
\coordinate (A)  at (0,5.2);
\coordinate (B)  at (-4.8,-2.8);
\coordinate (C)  at (4.8,-2.8);

\coordinate (D)  at (-2.4,1.2);
\coordinate (E)  at (2.4,1.2);
\coordinate (F)  at (0,-2.8);

\coordinate (AD) at (-1.2,3.2);
\coordinate (AE) at (1.2,3.2);
\coordinate (DE) at (0,1.2);

\coordinate (DB) at (-3.6,-0.8);
\coordinate (DF) at (-1.2,-0.8);
\coordinate (BF) at (-2.4,-2.8);

\coordinate (EF) at (1.2,-0.8);
\coordinate (EC) at (3.6,-0.8);
\coordinate (FC) at (2.4,-2.8);

% Resistor edges: upper cell
\draw (A)  to[R] (AD);
\draw (AD) to[R] (D);
\draw (A)  to[R] (AE);
\draw (AE) to[R] (E);
\draw (D)  to[R] (DE);
\draw (DE) to[R] (E);
\draw (AD) to[R] (AE);
\draw (AD) to[R] (DE);
\draw (AE) to[R] (DE);

% Resistor edges: lower-left cell
\draw (D)  to[R] (DB);
\draw (DB) to[R] (B);
\draw (D)  to[R] (DF);
\draw (DF) to[R] (F);
\draw (B)  to[R] (BF);
\draw (BF) to[R] (F);
\draw (DB) to[R] (DF);
\draw (DB) to[R] (BF);
\draw (DF) to[R] (BF);

% Resistor edges: lower-right cell
\draw (E)  to[R] (EF);
\draw (EF) to[R] (F);
\draw (E)  to[R] (EC);
\draw (EC) to[R] (C);
\draw (F)  to[R] (FC);
\draw (FC) to[R] (C);
\draw (EF) to[R] (EC);
\draw (EF) to[R] (FC);
\draw (EC) to[R] (FC);

% Junction nodes
\foreach \p in {A,B,C,D,E,F,AD,AE,DE,DB,DF,BF,EF,EC,FC}
  \fill (\p) circle (2.1pt);

% Current injection at B
\draw (B) -- ++(0,-0.55)
      to[I,l_=$i$] ++(0,-1.45)
      node[ground]{};
\node[blue!70!black, left=5pt] at (-4.8,-3.55) {$\text{source}$};

% Ground/reference at C
\draw (C) -- ++(0,-0.75) node[ground]{};
\node[green!50!black, right=7pt] at (4.8,-3.45) {$\text{ground / sink}$};

% Sampled voltage nodes
\fill[red!80!black] (AD) circle (4.0pt);
\node[red!80!black, left=5pt] at (AD) {$x_1$};

\fill[red!80!black] (AE) circle (4.0pt);
\node[red!80!black, right=5pt] at (AE) {$x_2$};

\fill[red!80!black] (DB) circle (4.0pt);
\node[red!80!black, left=5pt] at (DB) {$x_3$};

\fill[red!80!black] (EC) circle (4.0pt);
\node[red!80!black, right=5pt] at (EC) {$x_4$};

\fill[red!80!black] (BF) circle (4.0pt);
\node[red!80!black, below=5pt] at (BF) {$x_q$};

% Intended/perturbed sensor pair
\fill[red!80!black] (F) circle (4.2pt);
\node[red!80!black, above right=4pt] at (F) {$x_k$};

\fill[blue!75!black] (DF) circle (4.2pt);
\node[blue!75!black, below left=4pt] at (DF) {$y_k$};

\draw[dashed, gray!70, bend right=35] (F) to (DF);
\node[font=\large, fill=white, inner sep=1pt] at (0,0.25)
  {$R_{\mathrm{eff}}(x_k,y_k)$};

% Small legend
\begin{scope}[scale=1.75, shift={(6.25,1.5)}]
  \fill[red!80!black] (0,0) circle (4pt);
  \node[anchor=west] at (0.28,0) {sampled voltage node $x_j$};

  \fill[blue!75!black] (0,-0.65) circle (4pt);
  \node[anchor=west] at (0.28,-0.65) {perturbed location $y_k$};

  \fill (0,-1.30) circle (2.1pt);
  \node[anchor=west] at (0.28,-1.30) {unmeasured node};

  \draw (-0.15,-2.00) to[R] (0.90,-2.00);
  \node[anchor=west] at (1.15,-2.00) {resistive edge};
\end{scope}

\end{circuitikz}
\caption{A finite Sierpi\'nski gasket with resistive edges. Current enters at the source and leaves at the grounded node. The red nodes mark voltage measurements, and $(x_k,y_k)$ illustrates a displaced sensor. Its measurement error is controlled by $R_{\mathrm{eff}}(x_k,y_k)$. Circuit models on Sierpi\'nski gasket geometries also occur in antenna analysis \cite{ArrighettiDeCupisGerosa}.}
\label{fig:electrical-gasket-sampling}
\end{figure}

\subsection{Resistive Networks and Effective Resistance}
Let $G=(V,E,c)$ be a finite connected conductance network with $n=\left|V\right|\geq2$ and symmetric conductances $c_{xy}=c_{yx}>0$ on the edges. The physical weighted graph Laplacian is
\begin{equation}
(L_Gv)(x) \ = \ \sum_{y\sim x}c_{xy}\bigl(v(x)-v(y)\bigr),
\end{equation}
and its Dirichlet energy is
\begin{equation}
\mathcal E_G(v,v) \ = \ \sum_{\{x,y\}\in E}c_{xy}\left|v(x)-v(y)\right|^2 \ = \ v^{*}L_Gv.
\end{equation}
For a voltage potential $v$, this is exactly the dissipated power, which we write as $P_{\mathrm{diss}}(v)$. We place the uniform probability measure $\mu(x)=1/n$ on $V$. With this normalization the self-adjoint operator associated with $\mathcal E_G$ in $L^2(V,\mu)$ is not $L_G$ itself but
\begin{equation}
A \ = \ nL_G,
\end{equation}
since $\langle Au,v\rangle_{L^2(V,\mu)}=v^{*}L_Gu=\mathcal E_G(u,v)$. The corresponding resistance metric is the usual effective resistance
\begin{equation}
R_{\mathrm{eff}}(x,y) \ = \ (\mathbf e_x-\mathbf e_y)^{*}L_G^{+}(\mathbf e_x-\mathbf e_y),
\end{equation}
where $L_G^{+}$ denotes the Moore-Penrose pseudoinverse \cite{DoyleSnell,KleinRandic,LyonsPeres}.

This normalization identifies the geometric constant in the Fourier Ratio estimate with the classical Kirchhoff index.

\begin{proposition}\label{prop:electricaldictionary}
Let
\begin{equation}
0 \ = \ \nu_0<\nu_1\leq\cdots\leq\nu_{n-1}
\end{equation}
be the eigenvalues of $L_G$, and let $\lambda_j$ be the eigenvalues of the resistance-space operator $A=nL_G$. Then $\lambda_j=n\nu_j$ and
\begin{equation}
\sum_{j=1}^{n-1}\frac{1}{\lambda_j} \ = \ \frac{1}{n}\sum_{j=1}^{n-1}\frac{1}{\nu_j} \ = \ \frac{Kf(G)}{n^2} \ = \ \frac{1}{2n^2}\sum_{x,y\in V}R_{\mathrm{eff}}(x,y),
\end{equation}
where
\begin{equation}
Kf(G) \ = \ \sum_{\{x,y\}\subset V,\ x\neq y}R_{\mathrm{eff}}(x,y)
\end{equation}
is the Kirchhoff index of the network. Consequently the resistance constant from Theorem~\ref{thm:resistancegeometricrecovery} is
\begin{equation}
C_R \ = \ \frac{\sqrt{Kf(G)}}{n}.
\end{equation}
\end{proposition}

\begin{proof}
The identity $A=nL_G$ gives $\lambda_j=n\nu_j$. Since $\mu$ is uniform, Theorem~\ref{thm:resistancetrace} becomes
\begin{equation}
\sum_{j=1}^{n-1}\frac{1}{\lambda_j} \ = \ \frac{1}{2n^2}\sum_{x,y\in V}R_{\mathrm{eff}}(x,y).
\end{equation}
The ordered sum on the right counts every unordered pair twice, so it is $2Kf(G)$. This gives $\sum_{j=1}^{n-1}\lambda_j^{-1}=Kf(G)/n^2$, and substituting $\lambda_j=n\nu_j$ also gives the classical identity $Kf(G)=n\sum_{j=1}^{n-1}\nu_j^{-1}$ \cite{KleinRandic}. Taking square roots gives the formula for $C_R$.
\end{proof}

\subsection{Voltage Potentials and Current Injection}
We now write the Fourier Ratio estimate in electrical variables. Fix an orthonormal eigenbasis $\phi_0,\ldots,\phi_{n-1}$ of $L_G$ in the usual Euclidean inner product, with $\phi_0=n^{-1/2}\mathbf 1$. All Fourier Ratios in this section refer to this fixed choice. For a nonzero voltage potential $v$, write
\begin{equation}
\widetilde{v}(j) \ = \ \langle v,\phi_j\rangle_{\ell^2(V)}, \qquad \operatorname{FR}_G(v) \ = \ \frac{\sum_{j=0}^{n-1}\left|\widetilde{v}(j)\right|}{\left(\sum_{j=0}^{n-1}\left|\widetilde{v}(j)\right|^2\right)^{1/2}}.
\end{equation}
The $L^2(V,\mu)$ eigenfunctions are $e_j=\sqrt n\,\phi_j$, so the common factor relating the two coefficient sequences cancels in the ratio and $\operatorname{FR}_G$ is exactly the full-cutoff Fourier Ratio used above. Define
\begin{equation}
V_{\mathrm{mean}} \ = \ \frac1n\sum_{x\in V}v(x), \qquad V_{\mathrm{RMS}} \ = \ \left(\frac1n\sum_{x\in V}\left|v(x)\right|^2\right)^{1/2}.
\end{equation}

\begin{corollary}\label{cor:electricalFR}
Every nonzero voltage potential $v$ satisfies
\begin{equation}
\operatorname{FR}_G(v) \ \leq \ \frac{\left|V_{\mathrm{mean}}\right|}{V_{\mathrm{RMS}}}+\frac{\sqrt{Kf(G)}}{n}\frac{P_{\mathrm{diss}}(v)^{1/2}}{V_{\mathrm{RMS}}}.
\end{equation}
In the mean-zero voltage gauge, this simplifies to
\begin{equation}
\operatorname{FR}_G(v) \ \leq \ \frac{\sqrt{Kf(G)P_{\mathrm{diss}}(v)}}{nV_{\mathrm{RMS}}}.
\end{equation}
\end{corollary}

\begin{proof}
Apply Corollary~\ref{cor:uniformresistance} with $\left\|v\right\|_{L^2(V,\mu)}=V_{\mathrm{RMS}}$, $\int_V v\,d\mu=V_{\mathrm{mean}}$, and $\mathcal E_G(v,v)=P_{\mathrm{diss}}(v)$. Proposition~\ref{prop:electricaldictionary} identifies the remaining constant as $C_R=\sqrt{Kf(G)}/n$. Setting $V_{\mathrm{mean}}=0$ gives the second estimate.
\end{proof}

Voltage differences are unchanged when a constant is added to the potential. We choose the constant so that the mean voltage is zero. Let $i:V\rightarrow\mathbb{R}$ describe current injection at the nodes, with total current zero. Thus
\begin{equation}
\sum_{x\in V}i(x) \ = \ 0,
\end{equation}
and choose the mean-zero voltage gauge. The network equation is
\begin{equation}
L_Gv \ = \ i.
\end{equation}
Since $L_G$ is invertible on the mean-zero subspace, this determines
\begin{equation}
v \ = \ L_G^{+}i.
\end{equation}

\begin{proposition}\label{prop:electricalfilter}
Let $\widetilde{i}(j)=\langle i,\phi_j\rangle_{\ell^2(V)}$. For the mean-zero voltage potential generated by $i$,
\begin{equation}
v \ = \ \sum_{j=1}^{n-1}\frac{\widetilde{i}(j)}{\nu_j}\phi_j,
\end{equation}
and
\begin{equation}
P_{\mathrm{diss}}(v) \ = \ i^{*}L_G^{+}i \ = \ \sum_{j=1}^{n-1}\frac{\left|\widetilde{i}(j)\right|^2}{\nu_j}.
\end{equation}
Moreover,
\begin{equation}
V_{\mathrm{RMS}}^2 \ = \ \frac1n\sum_{j=1}^{n-1}\frac{\left|\widetilde{i}(j)\right|^2}{\nu_j^2}
\end{equation}
and, if $i\neq0$, the Fourier Ratio of the voltage has the exact current-source representation
\begin{equation}
\operatorname{FR}_G(v) \ = \ \frac{\sum_{j=1}^{n-1}\left|\widetilde{i}(j)\right|/\nu_j}{\left(\sum_{j=1}^{n-1}\left|\widetilde{i}(j)\right|^2/\nu_j^2\right)^{1/2}}.
\end{equation}
\end{proposition}

\begin{proof}
Expanding $i$ in the eigenbasis of $L_G$ and using $\widetilde{i}(0)=0$ gives
\begin{equation}
i \ = \ \sum_{j=1}^{n-1}\widetilde{i}(j)\phi_j.
\end{equation}
Applying $L_G^{+}$ divides the $j$th nonconstant coefficient by $\nu_j$, which gives the expansion of $v$. The power identity follows from
\begin{equation}
P_{\mathrm{diss}}(v) \ = \ v^{*}L_Gv \ = \ v^{*}i \ = \ i^{*}L_G^{+}i,
\end{equation}
and Parseval gives both spectral sums. The final formula is then the definition of $\operatorname{FR}_G(v)$ applied to the coefficients $\widetilde{i}(j)/\nu_j$.
\end{proof}

The network divides the current coefficient at frequency $j$ by $\nu_j$. This suppresses high frequencies in the voltage. It need not decrease the Fourier Ratio, since the ratio also involves an $\ell^2$ normalization. On a fixed spectral band we can compare the two ratios directly.

\begin{corollary}\label{cor:electricalFRdistortion}
Suppose the nonzero current injection $i$ has spectral support only on eigenvalues satisfying $0<\nu_a\leq\nu_j\leq\nu_b$. Then
\begin{equation}
\frac{\nu_a}{\nu_b}\operatorname{FR}_G(i) \ \leq \ \operatorname{FR}_G(v) \ \leq \ \frac{\nu_b}{\nu_a}\operatorname{FR}_G(i).
\end{equation}
\end{corollary}

\begin{proof}
On the spectral support of $i$ we have $\nu_b^{-1}\leq\nu_j^{-1}\leq\nu_a^{-1}$. Hence
\begin{equation}
\frac{1}{\nu_b}\left\|\widetilde{i}\right\|_1 \ \leq \ \left\|\widetilde{v}\right\|_1 \ \leq \ \frac{1}{\nu_a}\left\|\widetilde{i}\right\|_1
\end{equation}
and
\begin{equation}
\frac{1}{\nu_b}\left\|\widetilde{i}\right\|_2 \ \leq \ \left\|\widetilde{v}\right\|_2 \ \leq \ \frac{1}{\nu_a}\left\|\widetilde{i}\right\|_2.
\end{equation}
Divide the upper $\ell^1$ estimate by the lower $\ell^2$ estimate, and then divide the lower $\ell^1$ estimate by the upper $\ell^2$ estimate.
\end{proof}

The same spectral formula gives an approximation theorem for arbitrary current sources. The voltage need not belong to a prescribed spectral subspace, but its part above the cutoff can be estimated from the size of the current.

\begin{theorem}\label{thm:electricalsmoothing}
Let $i$ have total current zero, let $v=L_G^{+}i$ be the mean-zero voltage potential, and define
\begin{equation}
I_{\mathrm{RMS}} \ = \ \left(\frac1n\sum_{x\in V}\left|i(x)\right|^2\right)^{1/2}.
\end{equation}
For $\Lambda>0$, let $P_{\Lambda}$ denote the orthogonal projection onto the graph-spectral subspace spanned by the eigenvectors with $\nu_j\leq\Lambda$, and set $v_{\Lambda}=P_{\Lambda}v$. Then
\begin{equation}
\left\|v-v_{\Lambda}\right\|_{L^2(V,\mu)} \ \leq \ \frac{I_{\mathrm{RMS}}}{\Lambda},
\end{equation}
\begin{equation}
\mathcal E_G(v-v_{\Lambda},v-v_{\Lambda}) \ \leq \ \frac{nI_{\mathrm{RMS}}^2}{\Lambda},
\end{equation}
and, if $D_R=\max_{x,y}R_{\mathrm{eff}}(x,y)$,
\begin{equation}
\left\|v-v_{\Lambda}\right\|_{L^{\infty}(V)} \ \leq \ \left(\frac{nD_R}{\Lambda}\right)^{1/2}I_{\mathrm{RMS}}.
\end{equation}
In particular, if $i$ is graph-spectrally supported below $\Lambda$, then $v=v_{\Lambda}$.
\end{theorem}

\begin{proof}
By Proposition~\ref{prop:electricalfilter},
\begin{equation}
v-v_{\Lambda} \ = \ \sum_{\nu_j>\Lambda}\frac{\widetilde{i}(j)}{\nu_j}\phi_j.
\end{equation}
Therefore Parseval and $\nu_j>\Lambda$ give
\begin{equation}
\left\|v-v_{\Lambda}\right\|_{L^2(V,\mu)}^2 \ = \ \frac1n\sum_{\nu_j>\Lambda}\frac{\left|\widetilde{i}(j)\right|^2}{\nu_j^2} \ \leq \ \frac{1}{n\Lambda^2}\sum_{j=1}^{n-1}\left|\widetilde{i}(j)\right|^2 \ = \ \frac{I_{\mathrm{RMS}}^2}{\Lambda^2}.
\end{equation}
Similarly,
\begin{equation}
\mathcal E_G(v-v_{\Lambda},v-v_{\Lambda}) \ = \ \sum_{\nu_j>\Lambda}\frac{\left|\widetilde{i}(j)\right|^2}{\nu_j} \ \leq \ \frac{1}{\Lambda}\sum_{j=1}^{n-1}\left|\widetilde{i}(j)\right|^2 \ = \ \frac{nI_{\mathrm{RMS}}^2}{\Lambda}.
\end{equation}
The tail has mean zero, so for every $x\in V$ we may average its difference against the uniform measure and use the resistance inequality exactly as in Proposition~\ref{prop:resistancecoherence} to obtain
\begin{equation}
\left|(v-v_{\Lambda})(x)\right| \ \leq \ D_R^{1/2}\mathcal E_G(v-v_{\Lambda},v-v_{\Lambda})^{1/2} \ \leq \ \left(\frac{nD_R}{\Lambda}\right)^{1/2}I_{\mathrm{RMS}}.
\end{equation}
If $i$ has no spectral coefficients above $\Lambda$, then the expansion of Proposition~\ref{prop:electricalfilter} shows the same for $v$ and the tail vanishes.
\end{proof}

\subsection{Voltage Reconstruction and Sensor Stability}
We can now state the sampling theorem directly in electrical variables. For $\Lambda>0$, let
\begin{equation}
W_{\Lambda} \ = \ \operatorname{span}\{\phi_j:\nu_j\leq\Lambda\}, \qquad D_{\Lambda}^{G} \ = \ \#\{j:\nu_j\leq\Lambda\}.
\end{equation}
The corresponding cutoff for $A=nL_G$ is $n\Lambda$. Set
\begin{equation}
B_{\Lambda} \ = \ \max\left\{1,nD_R\Lambda\right\}.
\end{equation}
Proposition~\ref{prop:resistancecoherence} gives $K_{n\Lambda}^2\leq B_{\Lambda}$.

\begin{theorem}\label{thm:electricalrecovery}
Let $\Lambda\geq\nu_1$ and $0<\varepsilon<1/2$, let $v\in W_{\Lambda}$ be nonzero, and define
\begin{equation}
r_{\mathrm{el}}(v) \ = \ \frac{\left|V_{\mathrm{mean}}\right|}{V_{\mathrm{RMS}}}+\frac{\sqrt{Kf(G)}}{n}\frac{P_{\mathrm{diss}}(v)^{1/2}}{V_{\mathrm{RMS}}}.
\end{equation}
Choose $x_1,\ldots,x_q$ independently and uniformly from the nodes of $G$, and let $v^*\in W_{\Lambda}$ minimize the $\ell^1$ norm of its graph-spectral coefficients subject to
\begin{equation}
\left(\frac1q\sum_{k=1}^q\left|v^*(x_k)-v(x_k)\right|^2\right)^{1/2} \ \leq \ \varepsilon V_{\mathrm{RMS}}.
\end{equation}
If
\begin{equation}
q \ \geq \ C B_{\Lambda}\frac{r_{\mathrm{el}}(v)^2}{\varepsilon^2}\log^2\left(2+B_{\Lambda}\frac{r_{\mathrm{el}}(v)^2}{\varepsilon^2}\right)\log\left(2D_{\Lambda}^{G}\right),
\end{equation}
then with probability at least $1-(D_{\Lambda}^{G})^{-c}$,
\begin{equation}
\left\|v^*-v\right\|_{L^2(V,\mu)} \ \leq \ C_{\mathrm{rec}}\varepsilon V_{\mathrm{RMS}}.
\end{equation}
For mean-zero voltage potentials the first term in $r_{\mathrm{el}}(v)$ vanishes.
\end{theorem}

\begin{proof}
Apply Theorem~\ref{thm:resistancegeometricrecovery} to the finite measured resistance space associated with $G$. Proposition~\ref{prop:electricaldictionary} gives $C_R=\sqrt{Kf(G)}/n$, the cutoff $L$ from that theorem is $n\Lambda$, and $D_L=D_{\Lambda}^{G}$. The resistance-space complexity $r_R(v)$ is therefore exactly $r_{\mathrm{el}}(v)$, which gives the stated sample count and recovery estimate.
\end{proof}

Theorem~\ref{thm:electricalrecovery} expresses a sufficient sampling bound in electrical variables. Replacing $B_{\Lambda}$ by the actual bound
\begin{equation}
K_{n\Lambda}^{2} \ = \ n\max_{\nu_j\leq\Lambda}\max_{x\in V}\left|\phi_j(x)\right|^{2}
\end{equation}
in Theorem~\ref{thm:spectralrecovery} can give a smaller sample count. The same theorem also permits the actual Fourier Ratio in place of $r_{\mathrm{el}}(v)$. These distinctions matter when the goal is reconstruction from fewer measurements than $D_{\Lambda}^{G}$.

For a voltage generated by a current source, Theorem~\ref{thm:electricalsmoothing} controls the error made by first projecting to $W_{\Lambda}$. We may treat the remaining part of the voltage as measurement noise.

\begin{theorem}\label{thm:electricalapproximaterecovery}
Let $0<\varepsilon<1/2$, let $i$ have total current zero, let $v=L_G^{+}i$, and suppose $v_{\Lambda}=P_{\Lambda}v\neq0$. Define
\begin{equation}
\eta_{\Lambda} \ = \ \left(\frac{nD_R}{\Lambda}\right)^{1/2}I_{\mathrm{RMS}}, \qquad r_{\mathrm{el},\Lambda} \ = \ \frac{\sqrt{Kf(G)}}{n}\frac{P_{\mathrm{diss}}(v)^{1/2}}{\left\|v_{\Lambda}\right\|_{L^2(V,\mu)}}.
\end{equation}
Choose $x_1,\ldots,x_q$ independently and uniformly from $V$, and let $v^*\in W_{\Lambda}$ minimize the $\ell^1$ norm of its graph-spectral coefficients subject to
\begin{equation}
\left(\frac1q\sum_{k=1}^q\left|v^*(x_k)-v(x_k)\right|^2\right)^{1/2} \ \leq \ \eta_{\Lambda}.
\end{equation}
If
\begin{equation}
q \ \geq \ C B_{\Lambda}\frac{r_{\mathrm{el},\Lambda}^2}{\varepsilon^2}\log^2\left(2+B_{\Lambda}\frac{r_{\mathrm{el},\Lambda}^2}{\varepsilon^2}\right)\log\left(2D_{\Lambda}^{G}\right),
\end{equation}
then with probability at least $1-(D_{\Lambda}^{G})^{-c}$,
\begin{equation}
\left\|v^*-v\right\|_{L^2(V,\mu)} \ \leq \ C_{\mathrm{rec}}\varepsilon\left\|v_{\Lambda}\right\|_{L^2(V,\mu)}+C_{\mathrm{rec}}\eta_{\Lambda}+\frac{I_{\mathrm{RMS}}}{\Lambda}.
\end{equation}
\end{theorem}

\begin{proof}
Since $v_{\Lambda}$ is mean zero and $\mathcal E_G(v_{\Lambda},v_{\Lambda})\leq P_{\mathrm{diss}}(v)$, Corollary~\ref{cor:electricalFR} gives
\begin{equation}
\operatorname{FR}_G(v_{\Lambda}) \ \leq \ r_{\mathrm{el},\Lambda}.
\end{equation}
Theorem~\ref{thm:electricalsmoothing} gives
\begin{equation}
\left\|v-v_{\Lambda}\right\|_{L^{\infty}(V)} \ \leq \ \eta_{\Lambda},
\end{equation}
so $v_{\Lambda}$ is feasible for the minimization problem. Apply Theorem~\ref{thm:spectralrecovery} with target $v_{\Lambda}$, noise bound $\eta_{\Lambda}$, Fourier Ratio bound $r_{\mathrm{el},\Lambda}$, and $K_{n\Lambda}^2\leq B_{\Lambda}$. It follows that
\begin{equation}
\left\|v^*-v_{\Lambda}\right\|_{L^2(V,\mu)} \ \leq \ C_{\mathrm{rec}}\left(\varepsilon\left\|v_{\Lambda}\right\|_{L^2(V,\mu)}+\eta_{\Lambda}\right).
\end{equation}
Finally Theorem~\ref{thm:electricalsmoothing} also gives $\left\|v-v_{\Lambda}\right\|_{L^2(V,\mu)}\leq I_{\mathrm{RMS}}/\Lambda$, and the triangle inequality proves the result.
\end{proof}

Effective resistance also controls the effect of moving a voltage sensor from one node to another.

\begin{corollary}\label{cor:electricalsensor}
Let $v$ be a voltage potential and suppose intended sensor nodes $x_1,\ldots,x_q$ and actual sensor nodes $y_1,\ldots,y_q$ satisfy
\begin{equation}
R_{\mathrm{eff}}(x_k,y_k) \ \leq \ \delta
\end{equation}
for every $k$. Then
\begin{equation}
\left(\frac1q\sum_{k=1}^q\left|v(y_k)-v(x_k)\right|^2\right)^{1/2} \ \leq \ \delta^{1/2}P_{\mathrm{diss}}(v)^{1/2}.
\end{equation}
If $v\in W_{\Lambda}$, then the sufficient condition
\begin{equation}
\delta \ \leq \ \frac{\varepsilon^2}{n\Lambda}
\end{equation}
guarantees that this placement error is at most $\varepsilon V_{\mathrm{RMS}}$. For the reconstruction statement, assume $0<\varepsilon<1/2$ and $\Lambda\geq\nu_1$, choose the intended nodes independently and uniformly from $V$, and suppose $v\in W_{\Lambda}$ is nonzero. The observed values are $b_k=v(y_k)$. Put
\begin{equation}
\eta_{\mathrm{sens}} \ = \ \delta^{1/2}P_{\mathrm{diss}}(v)^{1/2},
\end{equation}
and let $v^*\in W_{\Lambda}$ minimize its spectral $\ell^1$ norm subject to
\begin{equation}
\left(\frac1q\sum_{k=1}^{q}\left|v^*(x_k)-b_k\right|^2\right)^{1/2} \ \leq \ \eta_{\mathrm{sens}}.
\end{equation}
Under the sample-size hypothesis of Theorem~\ref{thm:electricalrecovery}, with the same probability as in that theorem,
\begin{equation}
\left\|v^*-v\right\|_{L^2(V,\mu)} \ \leq \ C_{\mathrm{rec}}\left(\varepsilon V_{\mathrm{RMS}}+\delta^{1/2}P_{\mathrm{diss}}(v)^{1/2}\right).
\end{equation}
\end{corollary}

\begin{proof}
The first estimate is Theorem~\ref{thm:resistancesensor} with $R=R_{\mathrm{eff}}$ and $\mathcal E_G(v,v)=P_{\mathrm{diss}}(v)$. If $v\in W_{\Lambda}$, then
\begin{equation}
P_{\mathrm{diss}}(v) \ = \ \sum_{\nu_j\leq\Lambda}\nu_j\left|\widetilde{v}(j)\right|^2 \ \leq \ \Lambda\sum_{j=0}^{n-1}\left|\widetilde{v}(j)\right|^2 \ = \ n\Lambda V_{\mathrm{RMS}}^2,
\end{equation}
which gives the stated condition on $\delta$. The first estimate shows that $v$ is feasible for the stated minimization problem. Theorem~\ref{thm:spectralrecovery}, with measurement noise bounded by $\eta_{\mathrm{sens}}$, gives the final estimate under the normalization $A=nL_G$.
\end{proof}

The current-voltage equation and the resistance inequality play different roles. The first controls the voltage above a spectral cutoff through the eigenvalues of $L_G$. The second controls measurement errors through effective resistance and dissipated power. Combining them with spectral $\ell^1$ minimization yields the estimates above, with the limitations of the geometric sample count kept explicit.

\appendix
\section{Coefficient families and limiting distributions}
The observation in this appendix is independent of the resistance-space arguments. It applies when a coefficient vector is obtained by evaluating one continuous function along point sets with a limiting distribution.

For a nonzero vector $a\in\mathbb{C}^D$, write
\begin{equation}
\operatorname{FR}(a)=\frac{\left\|a\right\|_1}{\left\|a\right\|_2}.
\end{equation}
\begin{proposition}\label{prop:finiteFR}
Every nonzero $a\in\mathbb{C}^D$ satisfies
\begin{equation}
1\leq\operatorname{FR}(a)\leq\sqrt D,
\qquad
\operatorname{FR}(a)\geq\frac{\left\|a\right\|_2}{\left\|a\right\|_\infty}.
\end{equation}
\end{proposition}
\begin{proof}
The first bounds follow from $\left\|a\right\|_2\leq\left\|a\right\|_1$ and Cauchy--Schwarz. For the last bound, use
\begin{equation}
\left\|a\right\|_2^2\leq\left\|a\right\|_\infty\left\|a\right\|_1.
\end{equation}
\end{proof}

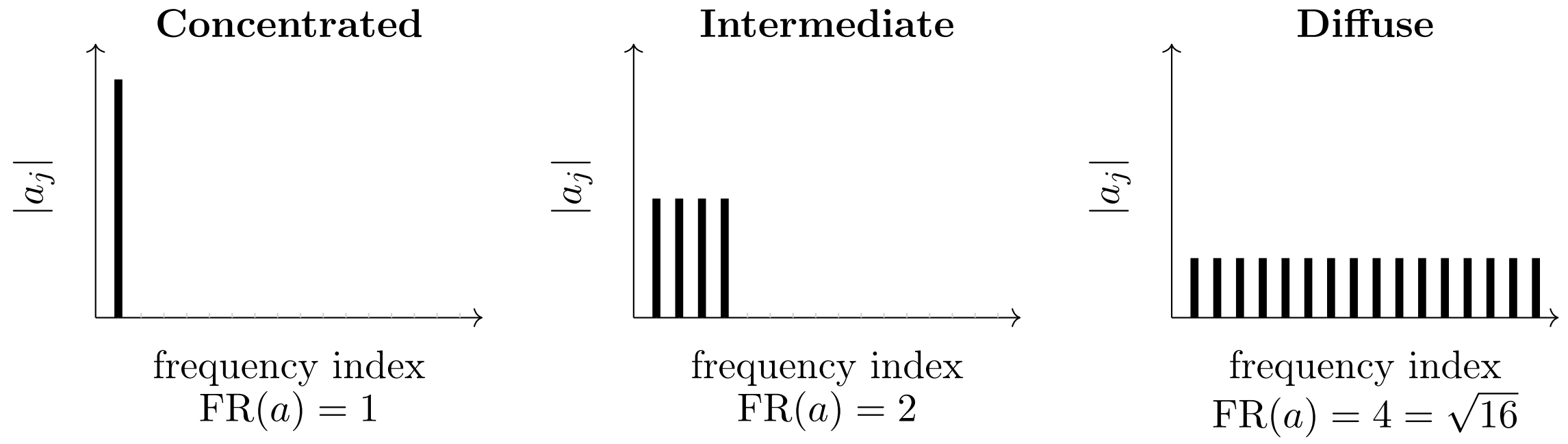
\begin{figure}[H]
\centering
\begin{tikzpicture}
\begin{scope}[xshift=0cm, x=0.22cm, y=2.3cm]
    \draw[->] (0,0) -- (17,0);
    \draw[->] (0,0) -- (0,1.15);

    \foreach \j in {1,...,16}{
        \draw[gray!35] (\j,0) -- (\j,0.02);
    }

    \draw[line width=2.2pt] (1,0) -- (1,1);

    \node[below] at (8.5,-0.08) {frequency index};
    \node[rotate=90, above] at (-1.2,0.55) {$|a_j|$};

    \node[above] at (8.5,1.12) {\textbf{Concentrated}};
    \node[below] at (8.5,-0.25)
        {$\operatorname{FR}(a)=1$};
\end{scope}

\begin{scope}[xshift=5.2cm, x=0.22cm, y=2.3cm]
    \draw[->] (0,0) -- (17,0);
    \draw[->] (0,0) -- (0,1.15);

    \foreach \j in {1,...,16}{
        \draw[gray!35] (\j,0) -- (\j,0.02);
    }

    \foreach \j in {1,...,4}{
        \draw[line width=2.2pt] (\j,0) -- (\j,0.5);
    }

    \node[below] at (8.5,-0.08) {frequency index};
    \node[rotate=90, above] at (-1.2,0.55) {$|a_j|$};

    \node[above] at (8.5,1.12) {\textbf{Intermediate}};
    \node[below] at (8.5,-0.25)
        {$\operatorname{FR}(a)=2$};
\end{scope}

\begin{scope}[xshift=10.4cm, x=0.22cm, y=2.3cm]
    \draw[->] (0,0) -- (17,0);
    \draw[->] (0,0) -- (0,1.15);

    \foreach \j in {1,...,16}{
        \draw[line width=2.2pt] (\j,0) -- (\j,0.25);
    }

    \node[below] at (8.5,-0.08) {frequency index};
    \node[rotate=90, above] at (-1.2,0.55) {$|a_j|$};

    \node[above] at (8.5,1.12) {\textbf{Diffuse}};
    \node[below] at (8.5,-0.25)
        {$\operatorname{FR}(a)=4=\sqrt{16}$};
\end{scope}
\end{tikzpicture}
\caption{The Fourier Ratio measures the effective spectral support of a coefficient vector. Each example has $\left\|a\right\|_2=1$, while the spectral mass is distributed over $1$, $4$, and $16$ equal coefficients, respectively. The corresponding Fourier Ratios are $1$, $2$, and $4$, illustrating the general range $1\leq\operatorname{FR}(a)\leq\sqrt D$.}
\label{fig:fourier-ratio-spectral-spread}
\end{figure}

\begin{theorem}\label{thm:equidistribution}
For every $D$, let $\beta_{D,1},\ldots,\beta_{D,D}$ lie in a fixed compact set $K\subset\mathbb{C}$. Suppose that the empirical probability measures
\begin{equation}
\nu_D=\frac1D\sum_{j=1}^D\delta_{\beta_{D,j}}
\end{equation}
converge weakly to a probability measure $\nu$ on $K$. Let $h:K\rightarrow\mathbb{C}$ be continuous and assume
\begin{equation}
\int_K\left|h(z)\right|^2\,d\nu(z)>0.
\end{equation}
Then
\begin{equation}
\frac{\sum_{j=1}^D\left|h(\beta_{D,j})\right|}
{\left(\sum_{j=1}^D\left|h(\beta_{D,j})\right|^2\right)^{1/2}}
=\sqrt D\,
\frac{\int_K\left|h(z)\right|\,d\nu(z)}
{\left(\int_K\left|h(z)\right|^2\,d\nu(z)\right)^{1/2}}
+o(\sqrt D).
\end{equation}
The left side is defined for all sufficiently large $D$.
\end{theorem}
\begin{proof}
Weak convergence applies to both continuous functions $|h|$ and $|h|^2$. Hence
\begin{equation}
\frac1D\sum_{j=1}^D\left|h(\beta_{D,j})\right|
\longrightarrow\int_K\left|h(z)\right|\,d\nu(z),
\end{equation}
and
\begin{equation}
\frac1D\sum_{j=1}^D\left|h(\beta_{D,j})\right|^2
\longrightarrow\int_K\left|h(z)\right|^2\,d\nu(z)>0.
\end{equation}
Divide the first average by the square root of the second and multiply by $\sqrt D$.
\end{proof}

\section*{Acknowledgments}
The authors gratefully acknowledge support from Williams College and the National Science Foundation under Grant No.~DMS-2241623 through the SMALL Undergraduate Research Program at Williams College.

\end{document}